\documentclass[noinfoline]{imsart}

\usepackage[english]{babel}
\RequirePackage[colorlinks]{hyperref} 

\usepackage{graphicx}
\usepackage{amsmath}
\usepackage{amsfonts}
\usepackage{amssymb}    
\usepackage{amsthm}
\usepackage{color}
\def\cov{\mbox{cov}}
\def\var{\mbox{Var}}

\usepackage{epsfig}

\def\1{{\mathbf 1}}

\newtheorem{thrm}{Theorem}[section]

\newtheorem{lemma}[thrm]{Lemma}

\newcommand{\pen}{\mathrm{pen}}

\begin{document}
\begin{frontmatter}
\title{Technical appendix to ``Adaptive estimation of stationary Gaussian fields''}
\runtitle{Technical appendix}

\begin{aug}
\author{\fnms{Nicolas} \snm{Verzelen}\ead[label=e3]{nicolas.verzelen@math.u-psud.fr}}
\address{Universit\'e Paris Sud, Laboratoire de Math\'ematiques, UMR 8628\\ Orsay Cedex F-91405\\ \printead{e3}}
\runauthor{Verzelen}
\end{aug}

\begin{abstract} 

This is a technical appendix to ``Adaptive estimation of stationary Gaussian fields''~\cite{verzelen_gmrf_theorie}. We present several proofs that have been skipped in the main paper. These proofs are organised as in Section 8 of \cite{verzelen_gmrf_theorie}.

\end{abstract}

 

\end{frontmatter}

\maketitle

\section{Proof of Proposition 8.1}

\begin{proof}[Proof of Proposition 8.1]

First, we recall the notations introduced in \cite{bousquet05}. Let $N$ be a positive integer.
Then, $\mathcal{I}_N$ stands for the family of subsets of $\{1,\ldots,N\}$
of size less than 2. Let $\mathcal{T}$  be a set of vectors indexed by
$\mathcal{I}_N$. In the sequel, $\mathcal{T}$ is assumed to be a compact subset of
$\mathbb{R}^{(N(N+1)/2)+1}$. The following lemma states
 a slightly modified version of the upper bound in remark $7$ in \cite{bousquet05}.

\begin{lemma}\label{lemmechaos_Rademacher}
Let $T$ be a supremum of Rademacher chaos indexed by $\mathcal{I}_N$ of the form
$$T:=\sup_{t\in\mathcal{T}}\bigg|\sum_{\{i,j\}}U_{i}U_{j}t_{\{i,j\}}+
  \sum_{i=1}^{N}t_{\{i\}} + t_{\varnothing}\bigg|\ ,$$
where $U_1,\ldots,U_{N}$ are independent Rademacher random variables.
Then for any $x>0$,
\begin{eqnarray}
\mathbb{P}\left\{T\geq \mathbb{E}[T]+x\right\} \leq
4\exp\left(-\frac{x^2}{L_1\mathbb{E}[D]^2}\wedge \frac{x}{L_2E}\right)\ ,
\end{eqnarray}
where $D$ and $E$ are defined by:
\begin{eqnarray}
D & := & \sup_{t\in \mathcal{T}}\sup_{\alpha:\|\alpha\|_2\leq
  1}\bigg|\sum_{i=1}^NU_i\sum_{j\neq i} \alpha_j t_{\{i,j\}}\bigg|\ , \nonumber\\
E & := & \sup_{t\in \mathcal{T}}\sup_{\alpha^{(1)},\alpha^{(2)},\|\alpha^{(1)}\|_2\leq 1\
  \|\alpha^{(2)}\|\leq 1}\bigg|\sum_{i=1}^N\sum_{j\neq i}t_{\{i,j\}}\alpha^{(1)}_i\alpha^{(2)}_j\bigg|\ . \nonumber
\end{eqnarray}
\end{lemma}
Contrary to the original result of \cite{bousquet05}, the chaos are not assumed to be homogeneous. Besides, the $t_{\{i\}}$ are redundant with $t_{\varnothing}$. In fact, we introduced this family in order to emphasize the connection with Gaussian chaos in the next result. 

A suitable application of the central limit theorem enables to obtain a corresponding bound for Gaussian chaos of order 2.
\begin{lemma}\label{lemmechaos1}
Let $T$ be a supremum of Gaussian chaos of order $2$.
\begin{eqnarray}
T:=\sup_{t\in \mathcal{T}}\bigg|\sum_{\{i,j\}} t_{\{i,j\}} Y_iY_j + \sum_{i}t_i Y_i^2 +t_{\varnothing}\bigg|\ ,
\end{eqnarray}
where $Y_1,\ldots, Y_N$ are independent standard Gaussian random variable. Then, for any $x>0$, 
\begin{eqnarray}
\mathbb{P}\left\{T\geq \mathbb{E}[T]+x\right\} \leq
\exp\left(-\frac{x^2}{\mathbb{E}[D]^2L_1}\wedge \frac{x}{EL_2}\right)\ ,
\end{eqnarray}
where	
\begin{eqnarray}
D & := & \sup_{t\in \mathcal{T}} \sup_{\alpha\in\mathbb{R}^N \|\alpha\|_2\leq
  1}\sum_{i,j} Y_i(1+\delta_{i,j})\alpha_jt_{\{i,j\}}\ ,
  \nonumber\\
E & := & \sup_{t\in \mathcal{T}}\sup_{\alpha_1,\ \|\alpha_1\|_2\leq
  1}\sup_{\alpha_2,\ \|\alpha_2\|_2\leq 1}\sum_{i,j}\alpha_{1,i}\alpha_{2,j}t_{\{i,j\}}(1+\delta_{i,j})\ .\nonumber
\end{eqnarray}
\end{lemma}
The proof of this Lemma is postponed to the end of this section. To conclude, we derive the result of Proposition 8.1 from this last lemma. For any matrix $R\in F$, we define the vector $t^{R}\in
\mathbb{R}^{nr(nr+1)/2+1}$ indexed by $\mathcal{I}_{nr}$
as follows
$$t_{\{(i,k),(j,l)\}}^{R}:= \delta_{k,l}(2-\delta_{i,j})\frac{R{\scriptstyle[i,j]}}{n},
\hspace{0.5cm} t_{\{(i,k)\}}^R := \frac{R{\scriptstyle[i,i]}}{n},\hspace{0.3cm}\text{ and }
t^R_{\varnothing}:=-tr(R)\ ,$$
where $\delta_{i,j}$ is the indicator function of $i=j$. In order to apply Lemma \ref{lemmechaos1} with $N=nr$ and
$\mathcal{T}=\left\{t^R|R\in F \right\}$,  we have to work out the quantities
$D$ and $E$.

\begin{eqnarray}
D & = & \sup_{t^R\in\mathcal{T}}\sup_{\alpha\in\mathbb{R}^{nr},\ \|\alpha\|_2\leq
  1} \left\{\sum_{i=1}^r\sum_{k=1}^n { Y}{\scriptstyle[i,k]} \sum_{j=1}^r\sum_{l=1}^n t^{R,k,l}_{ij}(1+\delta_{i,j}\delta_{k,l})\alpha_j^l  \right\} \nonumber\\
  & = & \sup_{R\in F}\sup_{\alpha\in\mathbb{R}^{nr},\ \|\alpha\|_2\leq 1}2
  \left\{\sum_{i=1}^r\sum_{k=1}^n{ Y}{\scriptstyle[i,k]} \sum_{j=1}^r \frac{R{\scriptstyle[i,j]}\alpha_j^k}{n} \right\} \nonumber \\ 
  & = & \sup_{R\in F}\sup_{\alpha\in\mathbb{R}^{nr},\ \|\alpha\|_2\leq
  1}\frac{2}{n} \left\{\sum_{k=1}^n\sum_{j=1}^r \alpha_j^k \left(\sum_{i=1}^r
  {Y}{\scriptstyle[i,k]} R{\scriptstyle[i,j]}\right)\right\}\ . \nonumber
\end{eqnarray}
Applying Cauchy-Schwarz identity yields 
\begin{eqnarray}
D^2& = & \frac{4}{n^2}\sup_{R\in F}\left\{\sum_{k=1}^n\sum_{j=1}^r
    \left(\sum_{i=1}^r {Y}{\scriptstyle [i,k]}R{\scriptstyle[i,j]}\right)^2\right\} \nonumber \\
    &  = & \frac{4}{n}\sup_{R\in F}tr(R\overline{YY^*}R^*)\ .\label{identite_D}
\end{eqnarray}
Let us now turn the constant $E$
\begin{eqnarray}
E & = & \sup_{t^R\in\mathcal{T}} \sup_{\begin{array}{c}\alpha_1,\alpha_2
 \in\mathbb{R}^{nr}\\ \|\alpha_1\|_2\leq 1,\|\alpha_2\|_2\leq 1\end{array}}\sum_{1\leq i,j\leq r}\sum_{1\leq k,l \leq n} (1+\delta_{ij}\delta_{k,l})t^{R,kl}_{i,j} \alpha^k_{1,i}\alpha^l_{2,j} \nonumber\\
 &  = & \sup_{R\in
 F}\sup_{\begin{array}{c}\alpha_1,\alpha_2\in\mathbb{R}^{nr}\\ \|\alpha_1\|_2\leq 1,\|\alpha_2\|_2\leq 1\end{array}}\frac{2}{n}\sum_{1\leq i,j\leq r}\sum_{1\leq k \leq n}R{\scriptstyle [i,j]}\alpha_{1,i}^k\alpha_{2,j}^k\ .\nonumber
\end{eqnarray}
From this last expression, it follows that $E$ is a supremum of $L_2$ operator norms
\begin{eqnarray*}
E= \frac{2}{n}\sup_{R\in F}\varphi_{\text{max}}\left( Diag^{(n)}(R)\right)\ ,
\end{eqnarray*}
where $Diag^{(n)}(R)$ is the $(nr\times nr)$ block diagonal
matrix such that each diagonal block is made of the matrix $R$. Since the largest eigenvalue of  $Diag^{(n)}(R)$ is exactly the largest eigenvalue of $R$, we get
\begin{eqnarray}\label{identite_E}
E= \frac{2}{n}\sup_{R\in F}\varphi_{\text{max}}(R)\ .
\end{eqnarray}   
Applying Proposition \ref{lemmechaos1} and gathering identities (\ref{identite_D}) and (\ref{identite_E}) yields
\begin{eqnarray*}
\mathbb{P}(Z \geq \mathbb{E}(Z) + t) \leq \exp\left[- \left( \frac{t^2}{L_1\mathbb{E}(V)}\bigwedge \frac{t}{L_2B}\right)\right]\ ,
\end{eqnarray*}
where $B=E$ and $V=D^2$.
\end{proof}\vspace*{0.5cm}

\begin{proof}[Proof of Lemma \ref{lemmechaos_Rademacher}]
This result is an extension of Corollary 4 in \cite{bousquet05}. We shall closely follow the sketch of
their proof adapting a few arguments. First, we upper bound the moments of
$\left(T-\mathbb{E}(T)\right)_+$. Then, we derive the deviation inequality from it. Here, $x_+=\max(x,0)$.
\begin{lemma}\label{majoration_momentq}
For all real numbers $q\geq 2$,  
\begin{eqnarray}\label{majo_momentq}
\|(T-\mathbb{E}(T))_+ \|_q\leq \sqrt{L q}\mathbb{E}(D) + L q E\ ,
\end{eqnarray}
where $\|T\|^q_q$ stands for the $q$-th moment of the random variable $T$. The quantities $D$ and $E$ are defined in Lemma \ref{lemmechaos_Rademacher}.
\end{lemma}
By Lemma \ref{majoration_momentq}, for any $t\geq 0 $ and any $q\geq 2$,
\begin{eqnarray*}
\mathbb{P}\left(T\geq \mathbb{E}(T) +t \right)& \leq&
\frac{\mathbb{E}\left[(T-\mathbb{E}(T))^q_+\right]}{t^q} \\
& \leq & \left(\frac{ \sqrt{L q}\mathbb{E}(D) + L q E}{t}\right)^q \ .
\end{eqnarray*}
The right-hand side is at most $2^{-q}$ if $\sqrt{L q}\mathbb{E}(D)\leq
t/4$ and $L qE\leq t/4$. Let us set 
$$q_0 := \frac{t^2}{16L \mathbb{E}(D)^2}\wedge \frac{t}{4L E}\ .$$
If $q_0\geq 2$, then $\mathbb{P}\left(T\geq \mathbb{E}(T) +t \right)\leq
2^{-q_0}$. On the other hand if $q_0<2$, then $4\times 2^{-q_0} \geq 1$. It follows that
\begin{eqnarray*}
\mathbb{P}\left(T\geq \mathbb{E}(T) +t \right)\leq 4 \exp\left(-\frac{\log(2)}{4L}\left[\frac{t^2}{4\mathbb{E}(D)^2}\wedge\frac{t}{E}\right]\right)\ .
\end{eqnarray*}
\end{proof}\vspace*{0.5cm}

\begin{proof}[Proof of Lemma \ref{majoration_momentq}]
This result is based on the entropy method developed in \cite{bousquet05}. Let $f:\mathbb{R}^N\rightarrow
\mathbb{R}$ be a  measurable function such that
$T=f(U_1,\ldots,U_N)$. In the sequel, $U'_1,\ldots, U'_N$ denote independent copies of
$U_1,\ldots,U_N$. The random variable $T'_i$ and $V^+$ are defined by
\begin{eqnarray*}
T'_i &:= &f(U_1,\ldots,U_{i-1},U'_i,U_{i+1},\ldots, U_N)\ ,\\
V^+ &:= &\mathbb{E}\left[\sum_{i=1}^N(T-T'_i)^2_+|U_1^N\right]\ ,
\end{eqnarray*}
where $U_1^N$ refers to the set $\left\{U_1,\ldots, U_N\right\}$.  Theorem 2 in \cite{bousquet05} states that for any real $q\geq 2$, 
\begin{eqnarray}\label{majorationmomentq}
\|(T-\mathbb{E}(T))_+ \|_q \leq \sqrt{Lq}\|\sqrt{V^+}\|_q\ .
\end{eqnarray}
To conclude, we only have bound the moments of $\sqrt{V^+}$. By definition,
$$T=\sup_{t\in\mathcal{T}}\bigg|\sum_{\{i,j\}}U_{i}U_{j}t_{\{i,j\}}+
  \sum_{i=1}^{N}t_{\{i\}} + t_{\varnothing}\bigg|\ .$$
Since the set $\mathcal{T}$ is compact, this supremum is achieved almost surely
  at an element $t^0$ of $\mathcal{T}$. For any $1\leq i\leq N$,
\begin{eqnarray*}
(T-T'_i)^2_+ \leq 
	 \bigg((U_i-U'_i)\bigg|\sum_{j\neq i}U_jt^0{\{i,j\}}\bigg|\bigg)^2\ .
\end{eqnarray*}
Gathering this bound for any $i$ between $1$ and $N$, we get
\begin{eqnarray*}
V^+ & \leq & \sum_{i=1}^N
\mathbb{E}\left[\left. \bigg((U_i-U'_i)\bigg|\sum_{j\neq
  i}U_jt^0{\{i,j\}}\bigg|\bigg)^2  \right|U_1^N\right]\\
& \leq & 2\sum_{i=1}^N \bigg[\sum_{j\neq
  i}U_jt^0{\{i,j\}}\bigg]^2 \\
& \leq & 2\sup_{\alpha\in \mathbb{R}^N,\ \|\alpha\|_2\leq
  1}\bigg[\sum_{i=1}^N\alpha_i\bigg(\sum_{j\neq
  i}t^0_{\{i,j\}}U_j\bigg)\bigg]^2 \\
& \leq & 2\sup_{t\in \mathcal{T}} \sup_{\alpha\in \mathbb{R}^N,\ \|\alpha\|_2\leq
  1}\sum_{i=1}^N \bigg[U_i\sum_{j\neq i}\alpha_jt_{\{i,j\}}\bigg]^2 =  2D^2\ .
\end{eqnarray*}
Combining this last bound with (\ref{majorationmomentq}) yields
\begin{eqnarray}
\|(T-\mathbb{E}(T))_+ \|_q & \leq & \sqrt{L q}\sqrt{2}\|D\|_q\nonumber\\
& \leq & \sqrt{L q}\left[\mathbb{E}(D) + \left|(D-\mathbb{E}(D))_+
\right\|_q\right]\ .\label{ineq_inter}
\end{eqnarray}

Since the random variable $D$  defined in Lemma \ref{lemmechaos_Rademacher} is a measurable function $f_2$ of the variables $U_1,\ldots, U_N$, we
apply again Theorem 2 in \cite{bousquet05}.
\begin{eqnarray*}
 \|(D- \mathbb{E}(D))_+\|_q \leq \sqrt{Lq}\left\|\sqrt{V^+_2}\right\|_q,
\end{eqnarray*}
where $V^+_2$ is defined by
\begin{eqnarray*}
V^+_2 := \mathbb{E}\left[\left.\sum_{i=1}^N(D-D'_i)_+^2\right|U_i^N\right],
\end{eqnarray*}
and $D'_i := f_2(U_1,\ldots,U_{i-1},U'_i,U_{i+1},\ldots, U_N)$. As previously, the supremum in $D$ is achieved at some
random parameter $(t^0,\alpha^0)$.  We therefore upper bound $V^+_2$ as previously.
\begin{eqnarray*}
V^+_2 & \leq & \sum_{i=1}^N\mathbb{E}\left[\bigg((U_i-U'_i)\bigg(\sum_{j\neq
  i}\alpha^0_jt^0_{\{i,j\}}\bigg)\bigg)^2\bigg|U_1^N\right]\\
& \leq &  2\sum_{i=1}^N\bigg(\sum_{j\neq i}\alpha^0_jt^0_{\{i,j\}}\bigg)^2\\
& \leq & 2\sup_{\alpha^{(2)}\in \mathbb{R}^N,\|\alpha\|_2\leq
  1}\bigg(\sum_{i=1}^N \alpha^{(2)}_j \sum_{j\neq i}\alpha^0_it_{\{i,j\}} \bigg)^2= 2 E^2\ .
\end{eqnarray*}
Gathering this upper bound with (\ref{ineq_inter}) yields 
$$\|(T-\mathbb{E}(T))_+ \|_q\leq \sqrt{L q}\mathbb{E}(D) + L q E\ .$$
\end{proof}

\begin{proof}[Proof of Lemma \ref{lemmechaos1}] We shall  apply the central limit theorem in order to transfer results for Rademacher chaos to Gaussian chaos.
Let $f$ be the unique function satisfying $T= f(y_1,\ldots ,y_N)$ for any $(y_1,\ldots,y_N)\in \mathbb{R}^N$. As the set $\mathcal{T}$ is compact, the function $f$ is known to be continuous. Let $(U^{(j)}_i)_{1\leq i\leq N,j\geq 0}$ an i.i.d. family of Rademacher
variables. For any integer $n>0$, the random variables $Y^{(n)}$ and $T^{(n)}$ are defined by
\begin{eqnarray}
Y^{(n)}& := &\bigg(\sum_{j=1}^n\frac{U_1^{(j)}}{\sqrt{n}},\ldots,
\sum_{j=1}^n\frac{U_N^{(j)}}{\sqrt{n}} \bigg	)\ ,\nonumber\\
T^{(n)}& := & f\left(Y^{(n)}\right)\ . \nonumber
\end{eqnarray}
Clearly, $T^{(n)}$ is a supremum of Rademacher chaos of order $2$ with $nN$
variables and a constant term.  By the central limit theorem, $T^{(n)}$ converges in distribution towards $T$ as $n$
tends to infinity. Consequently, deviation inequalities for the variables $T^{(n)}$ transfer to $T$ as long as the quantities $\mathbb{E}\left[D^{(n)}\right]$, $E^{(n)}$, and $\mathbb{E}[T^{(n)}]$ converge.

We first prove that the sequence $T^{(n)}$ converges in expectation towards $T$.
As $T^{(n)}$ converges in distribution, it is sufficient to show that the
sequence $T^{(n)}$ is asymptotically uniformly integrable. The set
$\mathcal{T}$ is compact, thus there exists a positive number $t_{\infty}$ such
that
\begin{eqnarray}
T^{(n)} & \leq &  t_{\infty}\bigg[\sum_{i,j}|Y_i^{(n)}Y_j^{(n)}|+1\bigg	]\nonumber \\
& \leq & t_{\infty}\bigg[1+(N+1)/2\sum_{i=1}^N \left(Y_i^{(n)}\right)^2 \bigg]\nonumber\ .
\end{eqnarray}
It follows that
\begin{eqnarray}
\left(T^{(n)}\right)^2 & \leq & t^2_{\infty}\left(\frac{N+1}{2}\right)^2\frac{N+2}{2}\left[1+\sum_{i=1}^{N}\left(Y_i^{(n)}\right)^4\right]\ . 
\end{eqnarray}

The sequence $Y_i^{(n)}$ does not only converge in distribution
to a standard normal distribution but also in moments (see for instance \cite{billingsley} p.391). It follows that 
$\overline{\lim} \mathbb{E}\left[
\left(T^{(n)}\right)^2\right] \leq \infty$
and the sequence $f\left(Y^{(n)}\right)$ is asymptotically uniformly integrable. As a consequence, $$\lim_{n\rightarrow \infty}\mathbb{E}\left[T^{(n)}\right] =\mathbb{E}[T]\ .$$

Let us turn to the limit of $\mathbb{E}\left[D^{(n)}\right]$. As the variable $T^{(n)}$ equals
\begin{eqnarray}
T^{(n)}=\sup_{t\in \mathcal{T}}\bigg|\sum_{\{i,j\}}t_{\{i,j\}}\sum_{1\leq k,l
  \leq n}\frac{U_i^{(k)}U_j^{(l)}}{n}+\sum_i t_i\sum_{1\leq k \leq n}
\frac{U_i^{(k)}}{\sqrt{n}}\sum_{l\neq
  k}\frac{U_i^{(l)}}{\sqrt{n}}+t_{\varnothing}+\sum_it_i\bigg|\ , \nonumber
\end{eqnarray}
it follows that 
\begin{eqnarray}
  D^{(n)} &  = &  \sup_{t\in\mathcal{T}}\sup_{\alpha\in\mathbb{R}^{nN},\ \|\alpha\|_2\leq 1}\bigg|\sum_{1\leq i\leq N}\sum_{1\leq k \leq n}U_i^{(k)}\bigg\{\sum_{j\neq i} \frac{t_{\{i,j\}}}{n}\sum_{1\leq l\leq n}\alpha_j^{(l)} + 2 \sum_{l\neq k} 2 \frac{t_{\{i\}}}{n}\alpha_i^{(l)} \bigg\}  \bigg| \nonumber\\
& \leq & \sup_{t\in\mathcal{T}}\sup_{\alpha\in\mathbb{R}^{nN},\ \|\alpha\|_2\leq
  1}\bigg\{\sum_i\frac{U_i^{(k)}}{\sqrt{n}}\sum_j
(1+\delta_{i,j})t_{\{i,j\}}\frac{\sum_{1\leq l\leq n}\alpha_j^{(l)}}{\sqrt{n}}
\bigg\} + A^{(n)}\ , \label{majoration_rn} 
\end{eqnarray}
where the random variable $A^{(n)}$ is defined by
$$A^{(n)} := \sup_{t\in\mathcal{T}}\sup_{\alpha\in\mathbb{R}^{nN},\ 
  \|\alpha\|_2\leq 1}\sum_{i=1}^N\sum_{j=1}^n t_{\{i\}} \frac{U^{(j)}_i}{n}\alpha^j_i\ .$$
Straightforwardly, one upper bounds $A^{(n)}$ by $t_{\infty}/n	\sqrt{\sum_{i=1}^N\sum_{j=1}^n\left(U_i^{(j)}\right)^2}$ and
  its expectation satisfies
$$\mathbb{E}\left(\left|A^{(n)}\right|\right)\leq t_{\infty} \sqrt{\frac{N}{n}}\ ,$$
which goes to 0 when $n$ goes to infinity. Thus, we only have to upper
bound the expectation of the first term in (\ref{majoration_rn}). Clearly, the supremum is achieved
only when for all $1\leq j\leq N$, the sequence 
$(\alpha_j^{(l)})_{1\leq l\leq n}$ is constant. In such a case, the sequence
$(\alpha_{j}^{(1)})_{1\leq j\leq N}$ satisfies  $\|\alpha^{(1)}\|_2 \leq
1/\sqrt{n}$. it follows that
\begin{eqnarray}
\mathbb{E}\left[D^{(n)}\right]= \mathbb{E}\bigg\{
\sup_{t\in\mathcal{T}}\sup_{\alpha\in\mathbb{R}^{N}\|\alpha\|_2\leq
  1}\mathbb{E}\bigg[\sum_i Y_i^{(n)}\sum_j(1+\delta_{i,j}) \alpha_j  \bigg]\bigg\} +
\mathcal{O}\left(\frac{1}{\sqrt{n}}\right)\ . \nonumber
\end{eqnarray} 
Let $g$ be the function defined by $$g\left(y_1,\ldots,y_N\right)=\sup_{t\in\mathcal{T}}\sup_{\alpha\in\mathbb{R}^{N}\|\alpha\|_2\leq
  1}\bigg[\sum_i y_i\sum_j(1+\delta_{i,j}) \alpha_j  \bigg]\ ,$$ for any $(y_1,\ldots,y_N)\in \mathbb{R}^N$. 
The function $g(.)$ is measurable and continuous as the supremum is taken over a compact set. As a consequence,
$g(Y^{(n)})$ converges in distribution towards $g(Y)$. As previously, the sequence is asymptotically
  uniformly integrable since its moment of order 2 is uniformly upper bounded. It follows that
$\lim \mathbb{E}\left[D^{(n)}\right] = \mathbb{E}\left[D\right]$.  ~\\

Third, we compute the limit of $E^{(n)}$. By definition, 
\begin{eqnarray}
E^{(n)} & = & \sup_{t\in\mathcal{T}} \sup_{\alpha_1,\alpha_2 \in
 \mathbb{R}^{nN},\ \|\alpha_1\|_2\leq 1,\|\alpha_2\|_2\leq 1}\sum_{i=1}^N\sum_{k=1}^n\alpha_{1,i}^k\bigg[\sum_{j\neq i}\sum_{l=1}^n \alpha_{2,j}^{(l)} \frac{t_{\{i,j\}}}{n}+ 2\sum_{l\neq k}\alpha_{2,i}^{(l)} \frac{t_{\{i\}}}{n} \bigg] \nonumber\\
 & = & \sup_{t\in\mathcal{T}} \sup_{\alpha_1,\alpha_2,\  \|\alpha_1\|_2\leq
 1,\|\alpha_2\|_2\leq
 1}\sum_{i=1}^N\sum_{j=1}^N(1+\delta_{i,j})\frac{t_{\{i,j\}}}{n}\left[\sum_{k=1}^n\sum_{l=1}^n\alpha_{1,i}^{(k)}\alpha_{2,j}^{(l)} \right] + \mathcal{O}\left(\frac{1}{n}\right)\ . \nonumber
\end{eqnarray}
As for the computation of $D^{(n)}$, the supremum is achieved when the sequences
$(\alpha_{1,i}^k)_{1\leq k\leq n}$ and $(\alpha_{2,j}^l)_{1\leq l\leq n}$ are constant for any $i\in
\{1,\ldots, N\}$. Thus, we only have to consider the supremum over the vectors
$\alpha_1$ and $\alpha_2$ in $\mathbb{R}^{N}$. 
\begin{eqnarray}
E^{(n)} = \sup_{t\in\mathcal{T}} \sup_{\alpha_1,\alpha_2
  \in\mathbb{R}^{N}\|\alpha_i\|_2\leq 1}\sum_{i=1}^N\sum_{j=1}^N
(1+\delta_{ij})t_{i,j} \alpha_{1,i}\alpha_{2,j} +\mathcal{O}\left(\frac{1}{n}\right). \nonumber
\end{eqnarray}
It follows that $E^{(n)}$ converges towards $E$ when $n$ tends to infinity.\vspace{0.5cm}~\\
The random variable $T^{(n)}-\mathbb{E}(T^{(n)})$ converges in distribution towards $T-\mathbb{E}(T)$. By Lemma \ref{lemmechaos_Rademacher} ,
\begin{eqnarray*}
\mathbb{P}(T-\mathbb{E}(T) \geq x ) \leq \underline{\lim}\  \exp\left(-\frac{x^2}{\mathbb{E}[D^{(n)}]^2L_1}\wedge \frac{x}{E^{(n)}L_2}\right)\ ,
\end{eqnarray*}
for any $x>0$. Combining this upper bound with the convergence of the sequences $D^{(n)}$ and $E^{(n)}$ allows to conclude.
\end{proof}\vspace*{0.5cm}

\section{Proof of Theorem 3.1}

\begin{proof}[Proof of Lemma 8.3]
We only consider here the anisotropic case, since the isotropic case is analogous.
This result is based on the deviation inequality for suprema of
Gaussian chaos of order $2$ stated in Proposition 8.1. For any model $m'$
belonging to $\mathcal{M}$, we shall upper bound the quantities
$\mathbb{E}(Z_{m'})$, $B_{m'}$, and $\mathbb{E}(W_{m'})$ defined in (42) in \cite{verzelen_gmrf_theorie}.

\begin{enumerate}
\item Let us first consider the expectation of $Z_{m'}$. Let $U'_{m,m'}$ be the new vector
space defined  by $$U'_{m,m'}:=U_{m,m'}\frac{\sqrt{D_{\Sigma}}}{p}\ , $$ where
$U_{m,m'}$ is introduced in the proof of Lemma 8.2 in \cite{verzelen_gmrf_theorie}. This new
space allows to handle the computation with the canonical inner product in the space of
matrices. Let $\mathcal{B}^{(2)}_{m^2,m'^2}$ be the unit ball of $U'_{m,m'}$ with respect
to the canonical inner product. If $R$ belongs to $U_{m,m'}$, then
$ \|R\|_{\mathcal{H}'}=\|R\sqrt{D_{\Sigma}}/p\|_F, $
where $\|.\|_F$ stands for the Frobenius norm. 
\begin{eqnarray}
Z_{m'}& =  & \sup_{R\in \mathcal{B}^{\mathcal{H}'}_{m^2,m'^2}}
 \frac{1}{p^2}tr\left[RD_{\Sigma}(\overline{\bf YY^*}-I_{p^2})\right] \nonumber\\
 & = & \sup_{R\in \mathcal{B}^{(2)}_{m^2,m'^2}}
 tr\left[R\frac{\sqrt{D_{\Sigma}}}{p}\left(\overline{\bf YY^*}-I_{p^2}\right)\right] \label{eqlemma}\\
 & =  & \bigg\|\Pi_{U'_{m,m'}}\frac{\sqrt{D_{\Sigma}}}{p}\left(\overline{\bf YY^*}-I_{p^2}\right)\bigg	\|_F\ , \nonumber
\end{eqnarray}  
where $\Pi_{U'_{m,m'}}$ refers to the orthogonal projection with respect to the canonical inner product onto the space $U'_{m,m'}$. 
Let $F_1,\ldots,F_{d_{m^2,m'^2}}$ denote an orthonormal basis of
$U'_{m,m'}$. 
\begin{eqnarray}
\mathbb{E}(Z_{m'}^2) & = & \sum_{i=1}^{d_{m^2,m'^2}}\mathbb{E}\bigg[
  tr^2\left(F_i\sqrt{\frac{D_{\Sigma}}{p^2}}\left(\overline{\bf YY^*}-I_{p^2}\right)\right)\bigg] \nonumber\\
& = & \sum_{i=1}^{d_{m^2,m'^2}}\mathbb{E}\bigg[
  \sum_{j=1}^{p^2}F_i{\scriptstyle[j,j]}\frac{\sqrt{D_{\Sigma}{\scriptstyle[j,j]}}}{p}(\overline{\bf YY^*}{\scriptstyle[j,j]}-1) \bigg]^2\nonumber\\	
 & =& \sum_{i=1}^{d_{m^2,m'^2}}\frac{2}{np^2}tr(F_iD_{\Sigma}F_i) \nonumber \\
 & \leq & \sum_{i=1}^{d_{m^2,m'^2}} \frac{2 \varphi_{\text{max}}\left( D_{\Sigma}\right) }{np^2}=
 \frac{2d_{m^2,m'^2}\varphi_{\text{max}}(\Sigma)}{np^2}\ .\nonumber 	
\end{eqnarray}
Applying Cauchy-Schwarz inequality, it follows that
\begin{eqnarray}
\mathbb{E}(Z_{m'})\leq \sqrt{ \frac{2d_{m^2,m'^2}\varphi_{\text{max}}(\Sigma)}{np^2}}\ .
\end{eqnarray}

\item Using the identity (\ref{eqlemma}), the quantity $B_{m'}$ equals
\begin{eqnarray}
B_{m'} &= & \frac{2}{n} \sup_{R\in \mathcal{B}^{(2)}_{m^2,m'^2}}\varphi_{\text{max}}\left( R\frac{\sqrt{D_{\Sigma}}}{p}  \right)\ . \nonumber
\end{eqnarray}
As the operator norm is under-multiplicative and as it dominates the Frobenius norm, we get the following bound
\begin{eqnarray}
 B_{m'} \leq   \frac{2\sqrt{\varphi_{\text{max}}(\Sigma)}}{np}\ .
\end{eqnarray}
\item Let us turn to bounding the quantity $\mathbb{E}(W_{m'})$. Again, by introducing the
 ball $\mathcal{B}^{(2)}_{m^2,m'^{2}}$, we get
\begin{eqnarray}
W_{m'} & = & \frac{4}{n}\sup_{R\in \mathcal{B}^{\mathcal{H}'}_{m^2,m'^2}}\frac{1}{p^2}tr\left[R\overline{\bf YY^*}D_{\Sigma}R\right] \nonumber \\
 & \leq & \frac{4\varphi_{\text{max}}(\Sigma)}{np^2} \sup_{R\in \mathcal{B}^{(2)}_{m^2,m'^2}}tr\left[R\overline{\bf YY^*}R\right] \nonumber\\
 & \leq &  \frac{4\varphi_{\text{max}}(\Sigma)}{np^2} \bigg( 1 +  \sup_{R\in \mathcal{B}^{(2)}_{m^2,m'^2}}tr\left[R\left(\overline{\bf YY^*}-I_{p^2}\right)R\right]\bigg)\ . \nonumber
\end{eqnarray}

Let $F_1,\ldots
F_{d_{m^2,m'^2}}$ an orthonormal basis of $U'_{m,m'}$ and let
$\lambda$ be a vector in $\mathbb{R}^{d_{m^2,m'^2}}$. We write $\|\lambda\|_2$ for its $L_2$ norm.
\begin{eqnarray*}
\lefteqn{\mathbb{E} \bigg(\sup_{R\in \mathcal{B}^{(2)}_{m^2,m'^2}}tr\left[R\left(\overline{\bf YY^*}-I_{p^2}\right)R\right]^2 \bigg) } & &\\& = & \mathbb{E}\bigg(\sup_{\|\lambda\|_2\leq 1} \sum_{i,j =1}^{d_{m^2,m'^2}} \lambda_i\lambda_j tr\left[ F_iF_j(\overline{\bf YY^*}/n-I_{p^2})\right]\bigg)^2 \nonumber\\
& \leq & \sum_{i,j=1}^{d_{m^2,m'^2}} \mathbb{E}\left(tr\left[F_iF_j({\bf YY^*}/n-I_{p^2})\right]^2\right)\ . \nonumber
\end{eqnarray*} 
The second inequality is a consequence of Cauchy-Schwarz inequality in
$\mathbb{R}^{(d_{m^2,m'^2})^2}$ since the $l_2$ norm of the vector
$(\lambda_i\lambda_j)_{1\leq i,j\leq d_{m^2,m'^2}}\in \mathbb{R}^{d_{m^2,m'^2}^2} $ is bounded by 1. Since the
matrices $F_i$ are diagonal, we  get
\begin{eqnarray*}
\mathbb{E} \bigg(\sup_{R\in \mathcal{B}^{(2)}_{m^2,m'^2}}tr\left[R({\bf YY^*}/n-I)R\right]^2 \bigg)\leq \frac{2}{n}  \sum_{i,j=1}^{d_{m^2,m'^2}} \|F_iF_j\|_2^2\ .
\end{eqnarray*}

It remains to bound the norm of the products $F_iF_j$ for any $i,j$ between
$1$ and $d_{m^2,m'^2}$.
\begin{eqnarray*}
\sum_{i,j=1}^{d_{m^2,m'^2}} \|F_iF_j\|_2^2 & = & \sum_{i,j=1}^{d_{m^2,m'^2}} \sum_{k=1}^{p^2}F_i{\scriptstyle[k,k]}^2F_j{\scriptstyle[k,k]}^2 = \sum_{k=1}^{p^2} \left(\sum_{i=1}^{d_{m^2,m'^2}}F_i{\scriptstyle[k,k]}^2  \right)^2\ .
\end{eqnarray*}
For any $k\in \{1,\ldots,p^2\}$,
$\sum_{i=1}^{d_{m^2,m'^2}}F_i{\scriptstyle[k,k]}^2\leq 1 $ since $(F_1,\ldots,F_{d_{m^2,m'^2}})$ form an orthonormal
family. Hence, we get
\begin{eqnarray*}
\sum_{i,j=1}^{d_{m^2,m'^2}} \|F_iF_j\|_2^2&\leq & \sum_{k=1}^{p^2}\sum_{i=1}^{d_{m^2,m'^2}}F_i{\scriptstyle[k,k]}^2=d_{m^2,m'^2}\nonumber \ .
\end{eqnarray*}
 All in all, we have proved that
\begin{eqnarray}
\mathbb{E}(W_{m'})\leq \frac{4\varphi_{\text{max}}(\Sigma)}{np^2}\left[1+\sqrt{\frac{2d_{m^2,m'^2}}{n}}\right].
\end{eqnarray}
\end{enumerate}
Gathering these three bounds and applying Proposition 8.1 allows to obtain the following deviation inequality:~\\
$$\begin{array}{l}
\mathbb{P}  \left(Z_{m'} \geq    \sqrt{\frac{2 \varphi_{\text{max}}(\Sigma)
  }{n}}\left\{\sqrt{1+\alpha/2}\sqrt{d_{m^2,m'^2}}+ \xi \right\}
\right)\\
  \leq    \exp\left\{- \left[\frac{\left[\left(\sqrt{1+\alpha/2}-1\right)\sqrt{d_{m^2,m'^2}}+\xi\right]^2}{2L_1\left(1+\sqrt{2d_{m^2,m'^2}/n}\right)} \bigwedge \frac{\sqrt{n}\left[\left(\sqrt{1+\alpha/2}-1\right)\sqrt{d_{m^2,m'^2}}+\xi\right]}{\sqrt{2}L_2}\right] \right\}\\
 \leq  \exp\left\{-
 \left[\frac{\omega_{m,m'}^2}{2L_1\left(1+\sqrt{2d_{m^2,m'^2}/n}\right)}
 \bigwedge
 \frac{\sqrt{n}\omega_{m,m'}}{\sqrt{2}L_2}\right]
 -
 \left[\frac{\xi\omega_{m,m'}}{L_1\left[1+\sqrt{2d_{m^2,m'^2}/n}\right]}
 \bigwedge \frac{\sqrt{n}\xi}{\sqrt{2}L_2}\right] \right\}\ ,
\end{array}$$ 
where $\omega_{m,m'}= \left(\sqrt{1+\alpha/2}-1\right)\sqrt{d_{m^2,m'^2}}$.
As $n$ and $d_{m^2,m'^2}$ are larger than one, there exists a
universal constant $L'_2$ such that
\begin{eqnarray*}
\left[\frac{(\sqrt{1+\alpha/2}-1)^2d_{m^2,m'^2}}{2L_1\left(1+\sqrt{2d_{m^2,m'^2}/n}\right)}
 \bigwedge
 \frac{\sqrt{n}(\sqrt{1+\alpha/2}-1)\sqrt{d_{m^2,m'^2}}}{\sqrt{2}L_2}\right]\hspace{3cm}\\ \hspace{3cm}\geq
 4L'_2\sqrt{d_{m^2,m'^2}}\left[\left(\sqrt{1+\alpha/2}-1\right)^2\wedge
 \left(\sqrt{1+\alpha/2}-1\right)\right]\ . 
\end{eqnarray*}
Since the vector space $U_{m,m'}$ contains all the matrices $D(\theta')$
with $\theta'$ belonging to $m'$, $d_{m^2,m'^2}$ is larger than
$d_{m'}$. Besides, by concavity of the square root function, it holds that $\sqrt{1+\alpha/2}-1\geq
\alpha[4\sqrt{1+\alpha/2}]^{-1}$. 
Setting 
$L'_1 := [4L_1(1+\sqrt{2})]^{-1} \wedge
[\sqrt{2}L_2]^{-1}$ and arguing as previously leads to
$$ \frac{\xi(\sqrt{1+\alpha/2}-1)\sqrt{d_{m^2,m'^2}}}{L_1\left(1+\sqrt{2d_{m^2,m'^2}/n}\right)}
 \bigwedge \frac{\sqrt{n}\xi}{\sqrt{2}L_2}\geq L'_1 \xi
 \left[\frac{\alpha}{\sqrt{1+\alpha/2}} \wedge \sqrt{n}\right]\ .$$
Gathering these two inequalities allows us to conclude that 
{\small
\begin{eqnarray*}
\lefteqn{\mathbb{P}  \left(Z_{m'} \geq    \sqrt{\frac{2 \varphi_{\text{max}}(\Sigma)
  }{n}}\left\{\sqrt{\left(1+\alpha/2\right)d_{m^2,m'^2}}+ \xi \right\}
\right)} &&\\&
 \leq &\exp\left\{-
 L'_2\sqrt{d_{m'}}\left(\frac{\alpha}{\sqrt{1+\alpha/2}}\wedge
 \frac{\alpha^2}{1+\alpha/2}\right) -
 L'_1\xi\left[\frac{\alpha}{\sqrt{1+\alpha/2}}\wedge \sqrt{n}\right] \right\}\ . 
\end{eqnarray*}}

\end{proof}\vspace*{0.5cm}

\begin{proof}[Proof of Lemma 8.4 in \cite{verzelen_gmrf_theorie}]
The approach falls in two parts. First, we relate the dimensions
$d_{m}$ and $d_{m^2}$ to the number of nodes of the torus $\Lambda$ that are closer
than $r_m$ or $2r_m$ to the origin $(0,0)$. We recall that the quantity $r_m$ is introduced in Definition 2.1 of \cite{verzelen_gmrf_theorie}. Second, we compute a nonasymptotic upper bound of the number of points in $\mathbb{Z}^2$ that lie in the disc of radius $r$. This second step is
quite tedious and will only give the main arguments.

Let $m$ be a model of the collection $\mathcal{M}_1$. By definition, $ m $ is the set of points lying in the disc of radius $r_m$ centered on $(0,0)$. Hence, 
\begin{eqnarray*}
\Theta_m = vect \left\{\Psi_{i,j},\ (i,j)\in  m   \right\}\ ,
\end{eqnarray*}
where the matrices $\Psi_{i,j}$ are defined by Eq. (14) in \cite{verzelen_gmrf_theorie}. As $\Psi_{i,j}=
\Psi_{-i,-j}$, the dimension $d_m$ of $\Theta_m$ is exactly the number of orbits of
$ m $ under the action of the central symmetry $s$. 

As $d_{m^2}$ is defined as the dimension of the space $U_m$, it also corresponds to
the dimension of the space
\begin{eqnarray}\label{espace2}
  vect\left\{C(\theta),\theta\in \Theta_m \right\}+
  vect\left\{C(\theta)^2,\theta\in \Theta_m\right\}\ ,
\end{eqnarray}
which is clearly in one to one correspondence with $U_m$. Straightforward computations lead to the following identity:
\begin{eqnarray}
C(\Psi_{i_1,j_1})C(\Psi_{i_2,j_2}) &= &C(\Psi_{i_1+i_2,j_1+j_2})\left[1+s_{i_1+i_2,j_1+j_2}\right]\nonumber\\ &+& C(\Psi_{i_1-i_2,j_1-j_2})\left[1+s_{i_1-i_2,j_1-j_2}\right])\ ,\nonumber
\end{eqnarray}
where $s_{x,y}$ is the indicator function of $x=-x$ and $y=-y$ in the torus $\Lambda$.
Combining this property with the definition of $\Theta_m$, we embed the space (\ref{espace2}) in the space
\begin{eqnarray*}
vect\left\{C\left(\Psi_{i_1+i_2,j_1+j_2}\right),\ (i_1,j_1),(i_2,j_2)\in  m \cup\{(0,0)\}\right\}\ ,
\end{eqnarray*}
and this last space is in one to one correspondence with 
\begin{eqnarray}
vect \left\{\Psi_{i_1+i_2,j_1+j_2},\  (i_1,j_1),(i_2,j_2)\in  m \cup\{(0,0)\} \right\}\ .
\end{eqnarray}
In the sequel, $\mathcal{N}(m)$ stands for the set $$\left\{(i_1+i_2,j_1+j_2),\  (i_1,j_1),(i_2,j_2)\in  m \cup\{(0,0)\}  \right\}\ .$$ Thus, the dimension $d_{m^2}$ is smaller or equal to the number of orbits of
$\mathcal{N}(m)$ under the action of the symmetry $s$.

To conclude, we have to compare the number of orbits in
$ m $ and the number of orbits in $\mathcal{N}(m)$. We distinguish two cases depending whether $2r_m+1\leq p$ or $2r_m+1> p$. First, we assume that $2r_m+1\leq p$. For such values the disc of radius $r_m$
centered on the points $(0,0)$ in not overlapping itself on the torus except
on a set of null Lebesgue measure.
\begin{figure}
\centerline{\epsfig{file=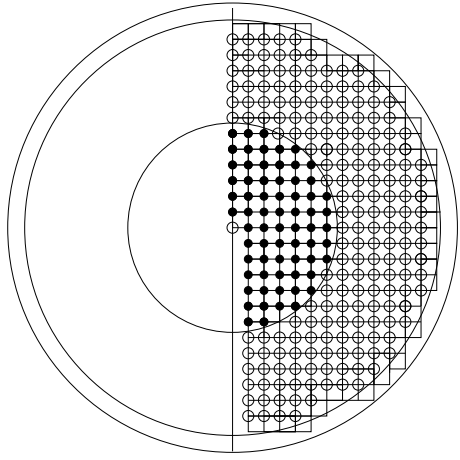,
    width=6cm}}
\caption{\textit{The black dots represent the orbit space of
    $ m $ and the white dots represent the remaining points of the orbit space
    of $\mathcal{N}(m)$.}\label{non-isotrope}}
\end{figure}
In the sequel, $\lfloor x\rfloor$ refers to the largest integer smaller than $x$.
We represent the orbit space of $ m $ as in Figure \ref{non-isotrope}. To any of
these points, we associate a square of size $1$. If we add $2+2\lfloor r_m\rfloor$ squares
to the $d_m$ first squares, we remark that the half disc centered on $(0,0)$ and
with length $r_m$ is contained in the reunion of these squares. Then, we get
\begin{eqnarray}
d_m+2 + 2\lfloor r_m\rfloor  \geq \frac{\pi r_m^2}{2}\ .  \label{majoration_dm1_anisotrope}
\end{eqnarray}

The points in $\mathcal{N}(m)$ are closer than $2r_m$ from the
origin. Consequently, all the squares associated to representants of
$\mathcal{N}(m)$ are included in the disc of radius $2r_m+\sqrt{2}$. 
\begin{eqnarray*}
d_{m^2} +2 +2 \lfloor 2r_m\rfloor \leq \frac{\pi}{2}\left\{ 2r_m+\sqrt{2}
\right\}^2\ .  \label{majoration_dm2_anisotrope}
\end{eqnarray*}
Combining these two inequalities, we are able to upper bound $d_{m^2}$ 
\begin{eqnarray}
2 +2\lfloor 2r_m \rfloor+d_{m^2}  &\leq & 4 \left\{1+\frac{\sqrt{2}}{2r_m} \right\}^2\left(d_m+1+2\lfloor r_m \rfloor\right) \nonumber\ ,\\
d_{m^2}& \leq & 4 \left\{1+\frac{\sqrt{2}}{2r_m} \right\}^2 d_m + 4
\left\{1+\frac{\sqrt{2}}{2r_m} \right\}^2\ (1+2\lfloor r_m \rfloor)\ .  \nonumber
\end{eqnarray}
Applying again inequality (\ref{majoration_dm1_anisotrope}), we upper bound
$r_m$:
\begin{eqnarray}
r_m \leq \frac{2}{\pi}\left[ 1+ \sqrt{1+\frac{\pi}{2}(1+d_m)}\right]\ . \nonumber
\end{eqnarray}
Gathering these two last bounds yields 
\begin{eqnarray}\label{majoration_ratio}
d_{m^2} \leq   4 \left\{1+\frac{\sqrt{2}}{2r_m} \right\}^2\left[1+
  \frac{1}{d_m}\left(1+\frac{4}{\pi}\left[ 1+
    \sqrt{1+\frac{\pi}{2}(1+d_m)}\right]\right)\right] d_m\ .  \nonumber
\end{eqnarray}

This upper bound is equivalent to $4d_m$, when $d_m$
goes to infinity. Computing the ratio $d_{m^2}/d_m$ for every model $m$
of small dimension allows to conclude.

Let us turn to the case $2r_m+1> p$. Suppose that $p$ is larger or equal to $9$. The lower bound
(\ref{majoration_dm1_anisotrope}) does not necessarily hold anymore. Indeed, the disc
is overlapping with itself because of toroidal effects. Nevertheless, we
obtain a similar lower bound by replacing $r_m$ by $(p-1)/2$:
\begin{eqnarray*}
d_m+2 + 2\lfloor \frac{p-1}{2}\rfloor  \geq \frac{\pi (p-1)^2}{8}.
\end{eqnarray*}
 
The number of orbits of $\Lambda$ under the action of the symmetry $s$ is
$(p^2+1)/2$ if $p$ is odd and $[(p+1)^2-1]/2$ if $p$ is
even. It follows that $d_{m^2}\leq [(p+1)^2-1]/2$.
Gathering these two bounds, we get
\begin{eqnarray*}
\frac{d_{m^2}}{d_m}\leq \frac{(p+1)^2}{\pi(p-1)^2/4-2(p+1)}\ .
\end{eqnarray*}
This last quantity is smaller than $4$ for any $p\geq 9$. An exhaustive
computation of the ratios when $p<9$ allows to conclude.\vspace{0.5cm}~\\

Let us turn to the isotropic case. Arguing as
previously, we observe that the dimension $d_m^{\text{iso}}$ is the number of orbits of the set
$ m $ under the action of the group $G$ introduced in
in \cite{verzelen_gmrf_theorie} Sect.1.1 whereas $d_{m^2}$ is smaller or equal to the number of orbits of
$\mathcal{N}^{\text{iso}}(m)$ under the action of $G$. As for anisotropic
models, we choose represent these orbits on the torus and associate squares of size $1$ (see Figure \ref{isotrope}).
\begin{figure}
\centerline{\epsfig{file=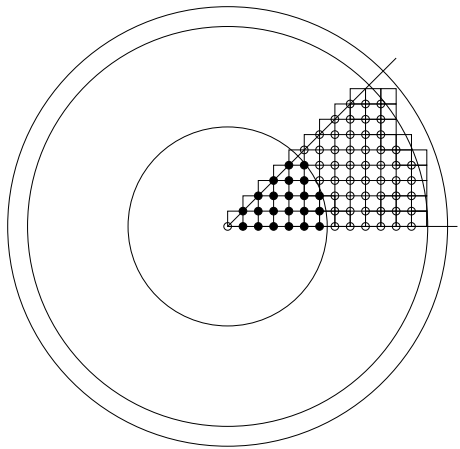,
    width=6cm}}
\caption{\textit{The black dots represent the orbit space of
    $ m $ under the action of $G$ and the white dots represent the remaining points of the orbit space
    of $\mathcal{N}^{\text{iso}}(m)$.
}\label{isotrope}}
\end{figure}
Assuming that $r_m<(p-1)/2$, we bound $d_m$ and $d_{m^2}$.

\begin{eqnarray*}\label{majoration_dm}
d_m +1  & \geq & \frac{1}{8}\pi r_m^2 + \frac{1}{2}\lfloor \frac{\sqrt{2}r_m}{2}\rfloor\ , \\
d_{m^2} & \leq &  4\left\{1+\frac{\sqrt{2}}{2r_m} \right\}^2\frac{1}{8}\pi r_m^2+\label{minoration_dmm}
\frac{1}{2}\lfloor \sqrt{2}r_m\rfloor\ .
\end{eqnarray*}
Gathering these two inequalities, we
get
$$ d_{m^2}\leq 4\left\{1+\frac{\sqrt{2}}{2r_m} \right\}^2 d_m\ . $$
As a consequence, $d_{m^2}$ is smaller than $4d_m$ when $d_m$ goes
to infinity. As previously, computing the ratio $d_{m^2}/d_m$ for models $m$ of small dimension allows to conclude.
The case $r_m> (p-1)/2$ is handled as for the anisotropic case. 
\end{proof}

\section{Proofs of the minimax bounds}

\begin{proof}[Proof of Lemma 8.5 in \cite{verzelen_gmrf_theorie}]
This lower bound is based on an application of Fano's approach. See \cite{yu} for a review of
this method and comparisons with Le Cam's and Assouad's Lemma. The proof follows three main steps:
First, we upper bound the Kullback-Leibler entropy between distributions
corresponding to $\theta_1$ and $\theta_2$ in the hypercube. Second, we find a set of points in the hypercube well
separated with respect to the Hamming distance. Finally, we conclude by applying Birg\'e's version
of Fano's lemma.  

\begin{lemma}\label{kullback}
The Kullback-Leibler entropy between two mean zero-Gaussian vectors of size $p^2$ with precision
matrices $\left(I_{p^2}-C(\theta_1)\right)/\sigma^2$ and $\left(I_{p^2}-C(\theta_2)\right)/\sigma^2$ equals
\begin{eqnarray}
 \mathcal{K}(\theta_1,\theta_2)= 1/2 \left[\log\left(\frac{|I_{p^2}-C(\theta_1)|}{|I_{p^2}-C(\theta_2)|}\right) +
 tr\left(\left[I_{p^2}-C(\theta_2)\right] \left[I_{p^2}-C(\theta_1)\right]^{-1}\right) - p^2 \right]\ , \nonumber
\end{eqnarray}
where for any square matrix $A$, $|A|$ refers to  the determinant of $A$.
\end{lemma}
This statement is classical and its proof is omitted. The
matrices $(I_{p^2}-C(\theta_1))$ and $(I_{p^2}-C(\theta_2))$ are diagonalizable in the same basis since they are symmetric block circulant (Lemma A.1 in \cite{verzelen_gmrf_theorie}). Transforming vectors of size $p^2$ into $p\times p$ matrices, we respectively define $\lambda_{1}$ and $\lambda_2$ as the $p\times p$ matrices of eigenvalues of $(I_{p^2}-C(\theta_1))$ and $(I_{p^2}-C(\theta_2))$. It follows that
$$\mathcal{K}(\theta_1,\theta_2)= 1/2 \sum_{1\leq i,j\leq p} \left(\frac{\lambda_2{\scriptstyle[i,j]}}{\lambda_1{\scriptstyle[i,j]}}-\log\left(\frac{\lambda_2{\scriptstyle[i,j]}}{\lambda_1{\scriptstyle[i,j]}}\right)-1\right)\ .$$
For any $x>0$, the following inequality holds
\begin{eqnarray}
x - 1 -\log(x) \leq \frac{9}{64}\left(x - \frac{1}{x}\right)^2. \nonumber
\end{eqnarray}
It is easy to establish by studying the derivative of corresponding functions.
As a consequence, 
\begin{eqnarray}
\frac{\lambda_2{\scriptstyle[i,j]}}{\lambda_1{\scriptstyle[i,j]}} -
 \log\left(\frac{\lambda_2{\scriptstyle[i,j]}}{\lambda_1{\scriptstyle[i,j]}}\right)- 1 
 & \leq &  \frac{9}{64}\left(\frac{\lambda_2{\scriptstyle[i,j]}}{\lambda_1{\scriptstyle[i,j]}}-
 \frac{\lambda_1{\scriptstyle[i,j]}}{\lambda_2{\scriptstyle[i,j]}}\right)^2 \nonumber \\ 
& \leq & \frac{9}{64}\left(\frac{1}{\lambda_1{\scriptstyle[i,j]}} + \frac{1}{\lambda_2{\scriptstyle[i,j]}}\right)^2\left(\lambda_1{\scriptstyle[i,j]}-\lambda_2{\scriptstyle[i,j]}\right)^2\ .
\end{eqnarray}
Let us first consider the anisotropic case. Let  $m$ be a model in $\mathcal{M}_1$ and let $\theta'$ belong $\Theta_m\cap\mathcal{B}_1(0_p,1)$. We also consider a positive radius $r$ such that $(1-\|\theta'\|_1-2rd_m)$ is positive.
For any $\theta_1$, $\theta_2$ in $\mathcal{C}_m(\theta',r)$ the matrices $(I_{p^2}-C(\theta_1))$ and $(I_{p^2}-C(\theta_2))$ are diagonally dominant and their eigenvalues  $\lambda_1{\scriptstyle[i,j]}$ and $\lambda_2{\scriptstyle[i,j]}$ are
larger than $1-\|\theta'\|_1-2rd_m$. 
\begin{eqnarray}
\mathcal{K}(\theta_1,\theta_2) & \leq
&\frac{9}{16(1-\|\theta'\|_1-2rd_m)^2}\sum_{1\leq i,j\leq p}(\lambda_{1}{\scriptstyle[i,j]}-\lambda_{2}{\scriptstyle[i,j]})^2
\nonumber \\ 
& \leq & \frac{9}{16(1-\|\theta'\|_1-2rd_m)^2} \| C(\theta_1)-C(\theta_2) \|_F^2
\nonumber \\
& \leq & \frac{9d_mr^2p^2}{8(1-\|\theta'\|_1-2rd_m)^2}\ . \label{majoration_kullback}
\end{eqnarray}
We recall that $\|.\|_F$ refers to the Frobenius norm in the space of
matrices.~\\

Let us state Birg\'e's version of Fano's lemma \cite{birgelemma}  and a combinatorial argument known under the name of Varshamov-Gilbert's lemma. These two lemma are taken from \cite{massartflour} and respectively correspond to Corollary 2.18 and Lemma 4.7.

\begin{lemma}\label{lemmebirge}\emph{\bf (Birg\'e's lemma)} Let $(S,d)$ be some pseudo-metric space and 
  $\{\mathbb{P}_s,s\in S\}$ be some statistical model. Let $\kappa$ denote some absolute constant smaller than one. Then for any estimator $\widehat{s}$ and any finite subset $T$ of $S$,
  setting $\delta=min_{s,t\in T,s\neq t}d(s,t)$, provided that
  $\max_{s,t\in T}\mathcal{K}(\mathbb{P}_s,\mathbb{P}_t)\leq \kappa \log|T|$, the
  following lower bound holds for every $p\geq 1$,
\begin{eqnarray*}
\sup_{s\in S}\mathbb{E}_s[d^p(s,\widehat{s})]\geq 2^{-p}\delta^p (1-\kappa)\ .
\end{eqnarray*}
\end{lemma}

\begin{lemma}\label{varshamov}\emph{\bf (Varshamov-Gilbert's lemma)}
Let $\{0,1 \}^d$ be equipped with Hamming distance $d_H$. There exists some
subset $\Phi$ of $\{0,1 \}^d$ with the following properties
\begin{eqnarray*}
  d_H(\phi,\phi')>d/4 \text{ for every $(\phi,\phi') \in \Phi^2$
  with $\phi\neq\phi'$}\text{ and }
\log|\Phi| \geq  \frac{d}{8}\ .
\end{eqnarray*}
\end{lemma}

Applying Lemma \ref{lemmebirge} with Hamming distance $d_H$ and the set $\Phi$ introduced 
in Lemma \ref{varshamov} yields
\begin{eqnarray}\label{minoration_hamming}
\sup_{\theta \in \mathcal{C}_m(\theta',r)} \mathbb{E}_{\theta}\left[d_H\left(\widehat{\theta},\theta\right)\right]\geq \frac{d_m}{8}(1-\kappa)\ ,
\end{eqnarray}
provided that
\begin{eqnarray}\label{condition_kulback}
\frac{9d_mr^2p^2n}{8(1-\|\theta'\|_1-2rd_m)^2} \leq \frac{\kappa d_m}{8}\ .
\end{eqnarray}
Let us express (\ref{minoration_hamming}) in terms of the Frobenius
 $\|.\|_F$ norm.
$$ \sup_{\theta \in \mathcal{C}_m(\theta',r)} \mathbb{E}_{\theta}\left[\|C(\widehat{\theta})-C(\theta)\|_F^2\right]\geq \frac{d_mr^2p^2}{4}(1-\kappa)\ . $$
Since for every $\theta$ in the hypercube, $\sigma^{-2}(I_{p^2}-C(\theta))$ is diagonally
dominant, its largest eigenvalue is smaller than $2\sigma^{-2}$. The loss function $l(\widehat{\theta},\theta)$ equals $\sigma^2/p^2tr\{[C(\widehat{\theta})-C(\theta)](I-C(\theta))^{-1}[C(\widehat{\theta})-C(\theta)]\}$. It follows that
\begin{eqnarray}\label{minoration_vrai}
\sup_{\theta \in \mathcal{C}_m(\theta',r)} \mathbb{E}_{\theta}\left[l\left(\widehat{\theta},\theta\right)\right]\geq \sigma^2\frac{d_mr^2}{8}(1-\kappa)\ .
\end{eqnarray} 
Condition (\ref{condition_kulback}) is equivalent to
$r^2(1-\|\theta'\|_1-2rd_m)^{-2} \leq \kappa /(9p^2n)$. If we assume that 
\begin{eqnarray}\label{condition3_kulback}
r^2\leq \frac{\kappa\left(1-\|\theta'\|_1\right)^2 }{18p^2n}\ , 
\end{eqnarray}
then $1-\|\theta'\|_1-2rd_m\geq \left(1-\|\theta'\|_1\right)\left(1-2d_m\sqrt{\kappa/(18np^2)}\right)$. This last quantity is
larger than $\left(1-\|\theta'\|_1\right)/\sqrt{2}$ if $d_m$ is smaller than $1.5(\sqrt{2}-1)\sqrt{np^2/\kappa}$. Gathering  inequality (\ref{minoration_vrai}) and condition (\ref{condition3_kulback}), we get the lower bound
\begin{eqnarray*}
\inf_{\widehat{\theta}}\sup_{\theta\in
  \text{Co}\left[\mathcal{C}_m(\theta',r)\right]}\mathbb{E}_{\theta}\left[l\left(\widehat{\theta},\theta\right)\right]
  &\geq& \inf_{\widehat{\theta}}\sup_{\theta\in\mathcal{C}_m\left[\theta',r\wedge
  (1-\|\theta'\|_1)\sqrt{\frac{\kappa}{18p^2n}}\right]}\mathbb{E}_{\theta}\left[l\left(\widehat{\theta},\theta\right)\right]\\ & \geq & L
  \left(r^2\wedge \frac{\left(1-\|\theta'\|_1\right)^2}{np^2}\right) d_m\sigma^2\ .
\end{eqnarray*}
One handles models of dimension $d_m$ between $1.5(\sqrt{2}-1)\sqrt{np^2/\kappa}$ and $\sqrt{n}p$ by changing the constant $L$ in the last lower bound.

\vspace{0.5cm}

Let us turn to sets of isotropic GMRFs. The proof is similar to the non-isotropic case, except for a few arguments.
Let $m$ belongs to the collection $\mathcal{M}_1$ and let $\theta'$ be an element of $\Theta_m^{\text{iso}}\cap\mathcal{B}_1(0_p,1)$. Let $r$ be such that $1-\|\theta'\|_1-8d_m^{\text{\emph{iso}}}$ is positive.
If $\theta_1$ and $\theta_2$ belong to the hypercube $\mathcal{C}^{\text{iso}}_{m}\left(\theta',r\right)$, then  \begin{eqnarray*}
\mathcal{K}(\theta_1,\theta_2) \leq \frac{9d_{m}r^2p^2}{2(1-\|\theta'\|_1-8rd_m^{\text{iso}})^2}\ .
\end{eqnarray*}
Applying Lemma \ref{lemmebirge} and \ref{varshamov}, it follows that
\begin{eqnarray*}
\inf_{\widehat{\theta}} \sup_{\theta \in C^{\text{iso}}_{m}(\theta',r)} \mathbb{E}_{\theta}\left[d_H\left(\widehat{\theta},\theta\right)\right]\geq \frac{d_{m}^{\text{iso}}}{8}(1-\kappa)\ ,
\end{eqnarray*}
provided that
$4.5d_{m}r^2p^2n(1-\|\theta'\|_1-8rd^{\text{iso}}_{m})^{-2}	 \leq \kappa d^{\text{iso}}_{m}/8$.
As a consequence, 
\begin{eqnarray*}
\inf_{\widehat{\theta}}\sup_{\theta \in C^{\text{iso}}_{m}(\theta',r)}
  \mathbb{E}_{\theta}\left[l\left(\widehat{\theta},\theta\right)\right]\geq \frac{d_{m}^{\text{iso}}r^2}{8}(1-\kappa)\ ,
\end{eqnarray*}
if $r^2\left(1-\|\theta'\|_1-8rd^{\text{iso}}_{m}\right)^{-2}\leq \kappa(36p^2n)^{-1}$. We conclude by arguing as in the isotropic case.
\end{proof}\vspace{0.5cm}

\begin{proof}[Proof of lemma 8.6 in \cite{verzelen_gmrf_theorie}]
Let $m$ be a model in $\mathcal{M}_1$, $r$ be a positive number smaller than  $1/(4d_m)$, and $\theta$ be an element of the convex hull of $\mathcal{C}_m(0_p,r)$. The covariance matrix of the vector $X^v$ is  $\Sigma= \sigma^2\left[I-C(\theta)\right]^{-1}$. Since the field $X$ is stationary, $\var_{\theta}(X{\scriptstyle[0,0]})$ equals any diagonal element of $\Sigma$. In particular, $\var_{\theta}(X{\scriptstyle[0,0]})$ corresponds to the mean of the eigenvalues of $\Sigma$. The matrix
$\left(I-C(\theta)\right)$ is  block circulant. As in the proof of Lemma \ref{minoration_hamming}, we note $\lambda$ the $p\times p$ matrix of the eigenvalues of ($I_{p^2}-C(\theta)$). By Lemma A.1 in \cite{verzelen_gmrf_theorie}, 
$$\lambda{\scriptstyle[i,j]}= 1 +
\sum_{(k,l)\in \Lambda}\theta{\scriptstyle[k,l]}\cos\left[2\pi\left(\frac{ik}{p}+\frac{jl}{p}\right)\right]\ ,$$
for any ${1\leq i,j\leq p}$. Since $\theta$ belongs to the convex hull of $\mathcal{C}_m(0_p,r)$,  $\theta{\scriptstyle[k,l]}$ is zero if
  $(k,l)\notin  m $ and $|\theta{\scriptstyle[k,l]}|\leq r$ if $(k,l)\in  m $. Thus $\sum_{(k,l)\in \Lambda}|\theta{\scriptstyle[k,l]}|$ is smaller than $1/2$. Applying Taylor-Lagrange inequality, we get 
$$\frac{1}{1+x}\leq 1 - x + \frac{x^2}{(1-|x|)^3}\ , $$
for any $x$ between $-1$ and $1$. It follows that
\begin{eqnarray*}
	\lambda{\scriptstyle[i,j]}^{-1} \leq 1 -	 
\sum_{k,l\in
  \Lambda}\theta{\scriptstyle[k,l]}\cos\left[2\pi\left(\frac{ik}{p}+\frac{jl}{p}\right)\right] +8\left\{\sum_{k,l\in \Lambda}\theta{\scriptstyle[k,l]}\cos\left[2\pi\left(\frac{ik}{p}+\frac{jl}{p}\right)\right]\right\}^2\ .
\end{eqnarray*}
Summing this inequality for all $(i,j)\in \left\{1,\ldots,p\right\}^2$, the first
order term turns out to be
$tr[C(\theta)]/p^2$ which is zero whereas the second term equals 
$8tr[C(\theta)^2]/p^2$. Since  there are less than $2d_m$ non-zero terms on each line of the matrix $C(\theta)$, its Frobenius norm is smaller than $2d_mp^2r^2$. Consequently, we obtain 
\begin{eqnarray*}
\var_{\theta}\left(X{\scriptstyle[0,0]}\right)\leq \sigma^2\left(1 + 16d_mr^2\right)\ .
\end{eqnarray*}

\end{proof}\vspace{0.5cm}

\begin{proof}[Proof of Lemma 8.7 in \cite{verzelen_gmrf_theorie}]
This property seems straightforward but the proof is a bit tedious.
Let $i$ be a positive integer smaller than $\text{Card}(\mathcal{M}_1)$.  By definition of the radius $r_m$ in  Equation (10) in \cite{verzelen_gmrf_theorie}, the model $m_{i+1}$ is the set of nodes in $\Lambda\setminus\{(0,0)\}$ at a distance smaller or equal to $r_{m_{i+1}}$ from $(0,0)$, whereas the model $m_i$ only contains the points in $\Lambda\setminus\{(0,0)\}$ at a distance strictly smaller than  $r_{m_{i+1}}$ from the origin. 

Let us first assume that $2r_{m_{i+1}}\leq p$. In such a case, the disc centered on $(0,0)$ with radius $r_{m_{i+1}}$ does not overlap with itself on the torus $\Lambda$. To any node in the neighborhood 
$m_{i+1}$ and to the node $(0,0)$, we
associate the square of size 1 centered on it. All these squares
do not overlap and are included in the disc of radius $r_{m_{i+1}}+\sqrt{2}/2$. Hence, we get the upper bound $
2d_{m_{i+1}}+1\leq \pi (r_{m_{i+1}}+\sqrt{2}/2)^2$. Similarly, the disc of radius $r_{m_{i+1}}-\sqrt{2}/2$ is included in the union of the squares associated to the nodes $m_i\cup\{0,0\}$. It follows that $2d_{m_i}+1$ is larger or equal to $\pi \left(r_{m_{i+1}}-\sqrt{2}/2\right)^2$. Gathering these two inequalities, we obtain
\begin{eqnarray*}
\frac{d_{m_{i+1}}}{d_{m_i}}\leq \frac{\left(r_{m_{i+1}}+\sqrt{2}/2\right)^2-1}{\left(r_{m_{i+1}}-\sqrt{2}/2\right)^2-1}\ ,
\end{eqnarray*}
if $r_{m_{i+1}}$ is larger than $1+\sqrt{2}/2$.
If $r_{m_{i+1}}$ larger than 5, this upper bound is smaller than two. An exhaustive computation for models of small dimension allows to conclude.

If $2r_{m_{i+1}}\geq p$ and  $2r_{m_i}< p$, then the preceding lower bound of $d_{m_i}$ and the preceding upper bound of $d_{m_{i+1}}$ still
hold. Finally, let us assume that $2r_{m_i}\geq p$. Arguing as previously, we conclude that
$2d_{m_{i}}+1\geq \pi (p/2 - \sqrt{2}/2)^2$. The largest dimension of a model $m\in \mathcal{M}_1$ is $(p^2-1)/2$ if $p$ is
odd and $((p+1)^2-3)/2$ if $p$ is even. Thus, $d_{m_{i+1}} \leq [(p+1)^2-3]/2$. Gathering these two bounds yields  
\begin{eqnarray*}
\frac{d_{m_{i+1}}}{d_{m_i}}\leq 4\frac{(p+1)^2-3}{\left(p-\sqrt{2}\right)^2}\ ,
\end{eqnarray*}
which is smaller than $2$ if $p$ is larger than $10$. Exhaustive computations for small $p$ allow to conclude.

\end{proof}\vspace{0.5cm}

\begin{proof}[Proof of Proposition 6.7 in \cite{verzelen_gmrf_theorie}]
This result derives from the upper bound of the risk of $\widetilde{\theta}_{\rho_1}$ stated in Theorem 3.1 and the minimax lower bound stated in Proposition 6.6 in \cite{verzelen_gmrf_theorie}. 

Let  $\mathcal{E}(a)$ be a pseudo-ellipsoid  that satisfies  Assumption $(\mathbb{H}_a)$ and such that $a_1^2\geq 1/(np^2)$.  For any $\theta$ in $\mathcal{E}(a)\cap\mathcal{B}_1(0_p,1)\cap\mathcal{U}(\rho_2)$, the penalty term satisfies $\pen(m)= K\sigma^2\rho_1^2\rho_2d_m/np^2$ is larger than  $Kd_m\varphi_{\text{max}}(\Sigma)/np^2$.
 Applying Theorem3.1, we upper bound the risk $\widetilde{\theta}_{\rho_1}$ 
\begin{eqnarray*}
\mathbb{E}_{\theta}\left[l\left(\widetilde{\theta}_{\rho_1},\theta\right)\right]\leq L_1(K)\inf_{m\in \mathcal{M}_1}
 \left[l(\theta_{m,\rho_1},\theta)+\pen(m)\right] +
L_2(K)\rho_2\frac{\sigma^2}{np^2}\ , 
\end{eqnarray*}
for any $\theta\in \mathcal{E}(a)\cap\mathcal{B}_1(0_p,1)\cap\mathcal{U}(\rho_2)$. It follows that
\begin{eqnarray*}
\sup_{\theta\in \mathcal{E}(a)\cap\mathcal{B}_1(0_p,1)\cap\mathcal{U}(\rho_2)}\mathbb{E}_{\theta}\left[l\left(\widetilde{\theta}_{\rho_1},\theta\right)\right]\leq L(K)\inf_{m\in \mathcal{M}_1\ ,\ d_m>0}
 \left[l(\theta_{m,\rho_1},\theta)+\rho_1^2 \rho_2\sigma^2\frac{d_m}{np^2}\right]\ .
\end{eqnarray*}

Let $i$ be a positive integer smaller or equal than  $\text{Card}(\mathcal{M}_1)$. We know from Section 4.1 in \cite{verzelen_gmrf_theorie} that the bias $l(\theta_{m_i},\theta)$ of the model $m_i$ equals $\var(X{\scriptstyle[0,0]}|X_{m_i}) -
  \sigma^2$. Since $\theta$ belongs to the set $\mathcal{E}(a)\cap\mathcal{B}_1(0_p,1)$, the bias term is smaller or equal to $a_{i+1}^2$ with the convention $a_{\text{Card}(\mathcal{M}_1)+1}^2= 0$. Hence , the previous upper bound becomes
\begin{eqnarray}
\mathbb{E}_{\theta}\left[l\left(\widetilde{\theta}_{\rho_1},\theta\right)\right]& \leq &  L(K)\inf_{1\leq i\leq
\text{Card}(\mathcal{M}_1)}\left[a_{i+1}^2 +\rho_1^2\rho_2\sigma^2\frac{d_{m_i}}{np^2}\right]\nonumber\\
& \leq &L(K,\rho_1,\rho_2)\inf_{1\leq i\leq\text{Card}(\mathcal{M}_1)}\left[a_{i+1}^2 +\frac{\sigma^2d_{m_i}}{np^2}\right]\label{majoration2_risque} \ . 
\end{eqnarray}
Applying Proposition 6.6 in \cite{verzelen_gmrf_theorie} to the set $\mathcal{E}(a)\cap\mathcal{B}_1(0_p,1)\cap\mathcal{U}(2)$, we get 
\begin{eqnarray*}
\inf_{\widehat{\theta}}\sup_{\theta\in\mathcal{E}(a)\cap\mathcal{B}_1(0_p,1)\cap\mathcal{U}(\rho_2)}\mathbb{E}_{\theta}\left[l\left(\widehat{\theta},\theta\right)\right]
& \geq &  \inf_{\widehat{\theta}}\sup_{\theta\in\mathcal{E}(a)\cap\mathcal{B}_1(0_p,1)\cap\mathcal{U}(2)}\mathbb{E}_{\theta}\left[l\left(\widehat{\theta},\theta\right)\right]\\
& \geq &  L
\sup_{1\leq i\leq \text{Card}(\mathcal{M}_1)}\left(a_i^2\wedge \sigma^2\frac{d_{m_i}}{np^2}\right) \ .
\end{eqnarray*}
Let us define $i^*$ by 
$$i^* :=  \sup\left\{1\leq i\leq \text{Card}(\mathcal{M}_1)\, ,\,a_i^2\geq \frac{\sigma^2d_{m_i}}{np^2}\right\}\ ,$$
with the convention $\sup \varnothing =0$. Since $a_1^2\geq \sigma^2/np^2$, $i^*$ is
larger or equal to one. It follows that 
\begin{eqnarray*}
\inf_{\widehat{\theta}}\sup_{\theta\in\mathcal{E}(a)\cap\mathcal{B}_1(0_p,\eta)}\mathbb{E}_{\theta}\left[l\left(\widehat{\theta},\theta\right)\right]\geq L_2
\left(a_{i^*+1}^2\vee \frac{\sigma^2d_{m_{i^*}}}{np^2}\right)\ .
\end{eqnarray*}
Meanwhile, the upper bound (\ref{majoration2_risque}) on the risk of $\widetilde{\theta}_{\rho_1}$  becomes
\begin{eqnarray*}
\mathbb{E}_{\theta}\left[l\left(\widetilde{\theta}_{\rho_1},\theta\right)\right]\leq L(K,\rho_1,\rho_2)\left( a_{i^*+1}^2 + \frac{\sigma^2d_{m_{i^*}}}{np^2}\right)
\leq  2L(K,\rho_1,\rho_2) \left(a_{i^*+1}^2\vee \frac{\sigma^2d_{m_{i^*}}}{np^2}\right)\ ,
\end{eqnarray*}
which allows to conclude.
\end{proof}

\section{Proof of the asymptotic risks bounds}

\begin{proof}[Proof of Corollary 4.6 in \cite{verzelen_gmrf_theorie}]

For the sake of simplicity, we assume that for any node $(i,j)\in m $, the nodes $(i,j)$ and $(-i,-j)$ are different in $\Lambda$. If this is not the case, we only have to slightly modify the proof in order to take account that $\|\Psi_{i,j}\|^2_F$ may equal one. The matrix $V$ is the covariance of the vector of size $d_m$ 
\begin{eqnarray}\label{cov_chi}
\left(X_{i_1,j_1}+X_{-i_1,-j_1},\ldots,X_{i_{d_m},j_{d_m}}+X_{-i_{d_m},-j_{d_m}} \right)\ .
\end{eqnarray}
Since the  matrix $\Sigma$ of $X^v$ is positive, $V$ is also positive. Moreover, its largest eigenvalue is larger than $2\varphi_{\text{max}}(\Sigma)$.

Let us assume first the $\theta$ belongs to $\Theta_m^+$ and that Assumption $(\mathbb{H}_1)$ is fulfilled. By the first result of Proposition 4.4 in \cite{verzelen_gmrf_theorie},
\begin{eqnarray*}
\lim_{n\rightarrow +\infty }np^2\mathbb{E}\left[l\left(\widehat{\theta}_{m,\rho_1},\theta\right)\right]=2\sigma^4tr\left[\text{\emph{IL}}_mV^{-1}\right]\geq \frac{\sigma^4}{\varphi_{\text{max}}(\Sigma)}tr[\text{\emph{IL}}_m]=2\sigma^4\frac{d_m}{\varphi_{\text{max}}(\Sigma)}\ ,
 \end{eqnarray*}
which corresponds to the first lower bound (30) in \cite{verzelen_gmrf_theorie}.\\

Let us turn to the second result. We now assume that $\theta$ satisfies Assumption $(\mathbb{H}_2)$. By the identity (28) of Proposition 4.4 in \cite{verzelen_gmrf_theorie}, we only have  to lower bound the quantity  $tr\left[VW^{-1}\right]$.
\begin{eqnarray*}
 tr\left[V^{-1}W\right]& \geq & \varphi_{\text{max}}(V)^{-1} tr\left[W\right]
 \geq  \frac{1}{2\varphi_{\text{max}}(\Sigma)}tr[W]\ .
\end{eqnarray*}
Since the matrix $\Sigma^{-1}=\sigma^{-2}\left[I_{p^2}-C(\theta)\right]$ is diagonally dominant, its smallest eigenvalue is  larger than $\sigma^{-2}(1-\|\theta\|_1)$. The matrix $\left(I_{p^2}-C(\theta_{m,\rho_1})\right)^2\left(I_{p^2}-C(\theta)\right)^{-2}$ is symmetric positive. It follows that $W$ is also symmetric positive definite. Hence, we get
\begin{eqnarray}
\lefteqn{tr\left[V^{-1}W\right]}\label{minoration_trace}& &\\& \geq &  \frac{\sigma^{-2}}{2}\left[1-\|\theta\|_1\right]\sum_{k=1}^{d_m}\frac{tr\left[C(\Psi_{i_k,j_k})^2\left[I_{p^2}-C(\theta_{m,\rho_1})\right]^2\left[I_{p^2}-C(\theta)\right]^{-2}\right]}{p^2}\nonumber . 
\end{eqnarray}
The 	largest eigenvalue of $\left(I_{p^2}-C(\theta)\right)$ is smaller than $2$ and the smallest eigenvalue of  $\left(I_{p^2}-C(\theta_{m,\rho_1})\right)$ is larger than $1-\|\theta_{m,\rho_1}\|_1$. By Lemma A.1 in \cite{verzelen_gmrf_theorie}, these two matrices are jointly diagonalizable and the smallest eigenvalue of $$\left(I_{p^2}-C(\theta_{m,\rho_1})\right)^2\left(I_{p^2}-C(\theta)\right)^{-2}$$	 is therefore larger than $(1-\|\theta_{m,\rho_1}\|_1)^2/4$. Gathering this lower bound with (\ref{minoration_trace}) yields 
\begin{eqnarray*}
 tr\left[V^{-1}W\right] \geq \frac{d_m\sigma^{-2}}{2}\left[1-\|\theta\|_1\right] \left[1-\|\theta_{m,\rho_1}\|_1\right]^2 \ .
\end{eqnarray*}
Lemma 4.1 in \cite{verzelen_gmrf_theorie} states that $\|\theta_{m,\rho_1}\|_1\leq \|\theta\|_1$. Combining these two lower bounds enables to conclude.

\end{proof}\vspace{0.5cm}

\begin{proof}[Proof of Example 4.8 in \cite{verzelen_gmrf_theorie}]

\begin{lemma}\label{variance_simple_lemme}
For any  $\theta$ is the space $\Theta_{m_1}^{+,\text{\emph{iso}}}$, the asymptotic variance term of $\widehat{\theta}_{m_1,\rho_1}^{\emph{\text{iso}}}$ equals 
\begin{eqnarray*} 
\lim_{n\rightarrow +\infty}np^2\mathbb{E}_{\theta}\left[l\left(\widehat{\theta}_{m_1,\rho_1}^{\emph{\text{iso}}},\theta \right)\right]= 2\sigma^4 \frac{tr\left(H^2\right)}{tr\left(H^2\Sigma\right)}\ .
\end{eqnarray*}
If $\theta$ belongs to $\Theta^{+,\text{\emph{iso}}}$ and also  satisfies $(\mathbb{H}_2)$, then 
\begin{eqnarray} \label{variance_asymptotique_simple2}
\lim_{n\rightarrow +\infty}np^2\mathbb{E}_{\theta}\left[l\left(\widehat{\theta}_{m_1,\rho_1}^{\emph{\text{iso}}},\theta_{m_1,\rho_1}^{\text{\emph{iso}}} \right)\right]= 2 \frac{tr\left\{\left[(I- \theta_{m_1,\rho_1}^{\text{\emph{iso}}}[1,0]H)H\Sigma\right]^2\right\}}{tr(H^2\Sigma)}\ ,
\end{eqnarray}
where the $p^2\times p^2$ matrix $H$ is defined as $H:=C\left(\Psi_{1,0}^{\text{\emph{iso}}}\right)$. 
\end{lemma}
\begin{proof}[Proof of Lemma \ref{variance_simple_lemme}]
Apply Proposition 4.4 in \cite{verzelen_gmrf_theorie} noting that $V= tr[H\Sigma H]/p^2$ and 
$$W =\frac{tr\left\{\left[(I- \theta_{m^{\text{iso}}_1}{\scriptstyle[1,0]}H)H\Sigma\right]^2\right\}}{\sigma^4p^2}\ . $$	
To prove the second result, we observe that $\Theta_{m_1}^{+,\text{iso}}$ equals $\Theta^{+,\text{iso}}_{m_1,2}$. It is stated for instance in Table 2 in \cite{verzelen_gmrf_theorie}.
\end{proof}

 Since the matrix $\theta$ belongs to $\Theta_{m_1}^{+,\text{iso}}$,  we may apply the second result of Lemma \ref{variance_simple_lemme}. Straightforward computations lead to $tr(H^2)=\|C\left(\Psi^{\text{iso}}_{1,0}\right)\|_F^2= 4p^2$ and $$tr(H^2\Sigma) = 4p^2\left[\var(X{\scriptstyle[0,0]})+ 2\cov_{\theta}(X{\scriptstyle[0,0]},X{\scriptstyle[1,1]})+ \cov_{\theta}\left(X{\scriptstyle[0,0]},X{\scriptstyle[2,0]}\right)\right] \ . $$
Since the field $X$ is an isotropic GMRF with four nearest neighbors, $$X{\scriptstyle[0,0]}= \theta{\scriptstyle[1,0]} \left(X{\scriptstyle[1,0]}+X{\scriptstyle[-1,0]}+ X{\scriptstyle[0,1]}+X{\scriptstyle[0,-1]}\right) + \epsilon{\scriptstyle[0,0]}\ , $$ where $\epsilon{\scriptstyle[0,0]}$ is independent from every variable $X{\scriptstyle[i,j]}$ with $(i,j)\neq 0$. Multiplying this identity by $X{\scriptstyle[1,0]}$ and taking the expectation yields
$$\cov_{\theta}\left(X{\scriptstyle[0,0]},X{\scriptstyle[1,0]}\right)= \theta{\scriptstyle[1,0]}\left[\var\left(X{\scriptstyle[0,0]}\right)+ 2 \cov_{\theta}(X{\scriptstyle[0,0]},X{\scriptstyle[1,1]})+ \cov_{\theta}\left(X{\scriptstyle[0,0]},X{\scriptstyle[2,0]}\right) \right]\ .$$
Hence, we obtain $tr\left(H^2\Sigma\right)= 4\cov_{\theta}(X{\scriptstyle[0,0]},X{\scriptstyle[1,0]})/\theta{\scriptstyle[1,0]}$ and 
$$\frac{tr(H^2)}{tr(H^2\Sigma)} = \frac{\theta{\scriptstyle[1,0]}}{\cov_{\theta}(X{\scriptstyle[0,0]},X{\scriptstyle[1,0]})}\ ,$$
which concludes the first part of the proof.\\

This second part is based on the spectral representation of the field $X$ and follows arguments which come back to Moran \cite{Moran73}.
We shall compute the limit of $\cov_{\theta}\left(X{\scriptstyle[0,0]},X{\scriptstyle[1,0]}\right)$ when the size of $\Lambda$ goes to infinity. As the field $X$ is stationary on $\Lambda$, we may diagonalize its covariance matrix $\Sigma$ applying Lemma A.1 in \cite{verzelen_gmrf_theorie}. We note $D_{\Sigma}$ the corresponding diagonal matrix defined by 
\begin{eqnarray*}
D_{\Sigma}{\scriptstyle[(i-1)p+j,(i-1)p+j]} = \sum_{k=1}^p\sum_{l=1}^p \cov_{\theta}\left(X{\scriptstyle[0,0]},X{\scriptstyle[k,l]}\right)\cos\left[2\pi\left(\frac{ki}{p}+\frac{lj}{p}\right)\right]\ ,
\end{eqnarray*}
for any $1\leq i,j\leq p$. Straightforwardly, we express $\cov_{\theta}\left(X{\scriptstyle[0,0]},X{\scriptstyle[1,0]}\right)$ as a linear combination of the eigenvalues
\begin{eqnarray*}
 \cov_{\theta}\left(X{\scriptstyle[0,0]},X{\scriptstyle[1,0]}\right) = \frac{1}{p^2} \sum_{i=1}^p\sum_{j=1}^p \cos\left(2\pi\frac{i}{p}\right)D_{\Sigma}{\scriptstyle[(i-1)p+j,(i-1)p+j]}\ .
\end{eqnarray*}
Applying Lemma A.1 in \cite{verzelen_gmrf_theorie} to the matrix $\Sigma^{-1}$ and noting that $\theta\in\Theta^{\text{iso},+}$ allows to get another expression of the eigenvalues of $\Sigma$
\begin{eqnarray*}
 D_{\Sigma}{\scriptstyle[(i-1)p+j,(i-1)p+j]} = \frac{\sigma^2}{1-2\theta{\scriptstyle[1,0]}\left[\cos\left(\frac{2\pi i}{p}\right)+ \cos\left(\frac{2\pi j}{p}\right)\right]}\ .
\end{eqnarray*}
We then combine these expression. By symmetry between $i$ and $j$ we get
\begin{eqnarray*}
 \cov_{\theta}\left(X{\scriptstyle[0,0]},X{\scriptstyle[1,0]}\right) = \frac{\sigma^2}{2p^2}\sum_{i=1}^p\sum_{j=1}^p \frac{\cos\left(2\pi\frac{i}{p}\right)+\cos\left(2\pi\frac{j}{p}\right)}{1-2\theta{\scriptstyle[1,0]}\left[\cos\left(2\pi \frac{i}{p}\right)+ \cos\left(2\pi \frac{j}{p}\right)\right]} \ .
\end{eqnarray*}
If we let $p$ go to infinity, this sum converges to the following integral
\begin{eqnarray*}
 \lefteqn{\lim_{p\rightarrow +\infty} \cov_{\theta}\left(X{\scriptstyle[0,0]},X{\scriptstyle[1,0]}\right)}&&\\& = & \frac{\sigma^2}{2}\int_{0}^1\int_0^1\frac{\cos(2\pi x)+ \cos(2\pi y)}{1-2\theta{\scriptstyle[1,0]}\left(\cos(2\pi x)+\cos(2\pi y)\right)}dxdy\\
& = & \frac{\sigma^2}{2\theta{\scriptstyle[1,0]}}\left[-1 + \frac{1}{4\pi^2}\int_0^{2\pi}\int_0^{2\pi}\frac{1}{1-2\theta{\scriptstyle[1,0]}\left[\cos(x)+\cos(y)\right]}dxdy\right]\ . 
\end{eqnarray*}
This last elliptic integral is asymptotically equivalent to $\log 16[4(1-4\theta{\scriptstyle[1,0]})]^{-1}$ when $\theta{\scriptstyle[1,0]}\rightarrow 1/4$ as observed for instance by Moran \cite{Moran73}. We conclude by substituting this limit in expression (33) in \cite{verzelen_gmrf_theorie}.
\end{proof}\vspace{0.5cm}

\begin{proof}[Proof of Example 4.9 in \cite{verzelen_gmrf_theorie}]
First, we compute $[\theta^{(p)}]_{m_1}^{\text{iso}}{\scriptstyle[1,0]}$. By Lemma 4.1 in \cite{verzelen_gmrf_theorie}, it minimizes
 the function $\gamma(.)$ defined in  (19) in \cite{verzelen_gmrf_theorie} over the whole space $\Theta_{m_1^{\text{iso}}}$. We therefore obtain
$$[\theta^{(p)}]_{m_1}^{\text{iso}}{\scriptstyle[1,0]} = \frac{tr\left[\Sigma H\right]}{tr\left[\Sigma H^2\right]}\ .$$
Once again, we apply Lemma A.1 in \cite{verzelen_gmrf_theorie} to simultaneously diagonalize the matrices $H$ and $\Sigma^{-1}$. As previously, we note $D_{\Sigma}$ the corresponding diagonal matrix of $\Sigma$. 
\begin{eqnarray*}
D_{\Sigma}{\scriptstyle[(i-1)p+j,(i-1)p+j]} & = & \frac{\sigma^2}{1-2\alpha\left[\cos\left(2\pi\left(\frac{pi}{4p}+\frac{pj}{4p}\right)\right)+\cos\left(2\pi\left(\frac{-pi}{4p}+\frac{pj}{4p}\right)\right)\right]}\\
& = & \frac{\sigma^2}{1-4\alpha\cos\left(\pi \frac{i}{2}\right)\cos\left(\pi\frac{j}{2}\right)} \ .
\end{eqnarray*}
Analogously, we compute the diagonal matrix $D\left(\Psi_{1,0}^{\text{iso}}\right)$
\begin{eqnarray*}
D\left(\Psi_{1,0}^{\text{iso}}\right){\scriptstyle[(i-1)p+j,(i-1)p+j]} & = 2\left[\cos\left(2\pi\frac{i}{p} \right) + \cos\left(2\pi\frac{j}{p} \right)\right] \ . 
\end{eqnarray*}
Combining these two last expressions, we obtain
\begin{eqnarray*}
 tr(H\Sigma) &= & \sum_{i=1}^{p}\sum_{j=1}^{p}\sigma^2 \frac{2\left[\cos\left(2\pi\frac{i}{p} \right) + \cos\left(2\pi\frac{j}{p} \right)\right]}{1-4\alpha\cos\left(\pi \frac{i}{2}\right)\cos\left(\pi\frac{j}{2}\right)} \ .
\end{eqnarray*}
Let us split this sum in 16 parts depending on the congruence of $i$ and $j$ modulo $4$. As each if of these 16 sums is shown to be zero, we conclude that  $tr(H\Sigma)= [\theta^{(p)}]_{m_1}^{\text{iso}}{\scriptstyle[1,0]}=0$.
By Lemma \ref{variance_simple_lemme}, the asymptotic risk of $\widehat{\theta^{(p)}}_{m_1}^{\text{iso},\rho_1}$ therefore equals
\begin{eqnarray*}
 \lim_{n\rightarrow +\infty} np^2 \mathbb{E}_{\theta^{(p)}}\left[l\left(\widehat{\theta^{(p)}}_{m_1}^{\text{iso},\rho_1},[\theta^{(p)}]_{m_1}^{\text{iso}}\right)\right] = \frac{tr(H^4\Sigma^2)}{tr(H^2\Sigma)}\ .
\end{eqnarray*}
First, we lower bound the numerator
\begin{eqnarray*}
 tr(H^4\Sigma^2) =\sigma^4 \sum_{i=1}^p\sum_{j=1}^p \frac{\left\{2\left[\cos\left(2\pi\frac{i}{p} \right) + \cos\left(2\pi\frac{j}{p} \right)\right]\right\}^4}{\left\{1-4\alpha\cos\left(\pi \frac{i}{2}\right)\cos\left(\pi\frac{j}{2}\right)\right\}^2}\ .
\end{eqnarray*}
As each term of this sum is non-negative, we may only consider the coefficients $i$ and $j$ which are congruent to 0 modulo $4$. \begin{eqnarray*}
 tr(H^4\Sigma^2) \geq\sigma^4 \sum_{i=0}^{p/4-1}\sum_{j=0}^{p/4-1}\frac{16 \left[\cos\left(2\pi\frac{i}{p/4} \right) + \cos\left(2\pi\frac{j}{p/4} \right)\right]^4}{(1-4\alpha)^2}\ .
\end{eqnarray*}
If we let go $p$ to infinity, we get the lower bound
\begin{eqnarray*}
\lim_{p\rightarrow +\infty}\frac{tr(H^4\Sigma^2)}{p^2}\geq \frac{\sigma^4}{(1-4\alpha)^2}\int_{0}^1\int_0^1 \left[\cos(2\pi x)+ \cos(2\pi y)\right]^4dxdy\ .
\end{eqnarray*}
Similarly, we upper bound $tr(H^2\Sigma)$ and let $p$ go to infinity
\begin{eqnarray*}
 \lim_{p\rightarrow +\infty}\frac{tr(H^2\Sigma)}{p^2}\leq \frac{4\sigma^2}{1-4\alpha}\int_{0}^1\int_0^1 \left[\cos(2\pi x)+ \cos(2\pi y)\right]^2 dxdy \ .
\end{eqnarray*}
Combining these two bounds allows to conclude 
\begin{eqnarray*}
\lim_{p\rightarrow +\infty}\lim_{n\rightarrow +\infty} np^2 R_{\theta^{(p)}}\left(\widehat{\theta^{(p)}}_{m_1}^{\text{iso},\rho_1},[\theta^{(p)}]_{m_1}^{\text{iso}}\right) \geq \frac{L\sigma^2}{1-4\alpha}\ .
\end{eqnarray*}
\end{proof}

\section{Miscellaneous}\label{appendix}

\begin{proof}[Proof of Lemma 1.1 in \cite{verzelen_gmrf_theorie}]
Let $\theta$ be a $p\times p$ matrix that satisfies condition
(3) in \cite{verzelen_gmrf_theorie}. For any $1\leq i_1,i_2\leq p$,
we define the $p\times p$ submatrix $C_{i_1,i_2}$ as
$$ C_{i_1,i_2}{\scriptstyle[j_1,j_2]} := C(\theta){\scriptstyle[(i_1-1)p+j_1,(i_2-1)p+j_2]}\ ,$$
for any $1\leq j_1,j_2\leq p$. For the sake of simplicity, the subscripts $(i_1,i_2)$ are taken modulo $p$. By definition of $C(\theta)$, it holds that $C_{i_1,i_2}=
C_{0 ,i_2-i_1}$ for any $1\leq i_1,i_2\leq p$. Besides, the matrices $C_{0 ,i}$ are circulant for any $1\leq i\leq p$. In short,
the matrix $C(\theta)$ is of the form
\begin{eqnarray*}
C(\theta) = \left(\begin{array}{cccc} 
C_{0,1} & C_{0,2} & \cdots & C_{0,p} \\
\vdots & \vdots & \vdots & \vdots\\
C_{0,p} & C_{0,1} & \cdots & C_{0,p-1}\end{array}
\right),
\end{eqnarray*}
where the matrices $C_{0,i}$ are circulant. Let $(i_1,i_2,j_1,j_2)$ be in $\{1,\ldots,p\}^4$. By definition, 
$$C(\theta){\scriptstyle[(i_1-1)p+j_1,(i_2-1)p+j_2]} = \theta{\scriptstyle[i_2-i_1,j_2-j_1]}\ .$$
Since the matrix $\theta$ satisfies condition (3) in \cite{verzelen_gmrf_theorie}, $\theta{\scriptstyle[i_2-i_1,j_2-j_1]}=
\theta{\scriptstyle[i_1-i_2,j_1-j_2]}$. As a consequence,~\\
$C(\theta){\scriptstyle[(i_1-1)p+j_1,(i_2-1)p+j_2]} = C(\theta){\scriptstyle[(i_2-1)p+j_2,(i_1-1)p+j_1]}$
and $C(\theta)$ is symmetric. \vspace{0.3cm}

Conversely, let $B$ be a $p^2\times p^2$ symmetric block circulant matrix. Let us define the matrix $\theta$ of size $p$ by
\begin{eqnarray*}
\theta{\scriptstyle[i,j]} := B{\scriptstyle[1,(i-1)p+j]}\ ,
\end{eqnarray*}
for any $1\leq i,j\leq p$.
Since the matrix $B$ is block circulant, it follows  that
$C(\theta)=B$. By definition,  $\theta{\scriptstyle[i,j]} = C(\theta){\scriptstyle[1,(i-1)p+j]}$ and $\theta{\scriptstyle[-i,-j]} = C(\theta){\scriptstyle[(i-1)p+j, 1 ]}$ for any integers $1\leq i,j\leq p$. Since the
matrix $B$ is symmetric, we conclude that $\theta{\scriptstyle[i,j]}=\theta{\scriptstyle[-i,-j]}$.
\end{proof}\vspace{0.5cm}

\begin{proof}[Proof of Lemma 2.2 in \cite{verzelen_gmrf_theorie}]
For any $\theta'\in \Theta^+$, $\gamma_{n,p}(\theta')$ is defined as
$$\gamma_{n,p}(\theta')=\frac{1}{p^2}tr\left[(I_{p^2}-C(\theta'))\overline{{\bf X^v X^{v*}}}(I_{p^2}-C(\theta'))\right]\ .$$
Applying Lemma A.1 in \cite{verzelen_gmrf_theorie}, there exists an orthogonal matrix $P$ that simultaneously diagonalizes $\Sigma$ and any matrix $C(\theta')$. Let us define $ {\bf Y}^i:= \sqrt{\Sigma}^{-1}{\bf X}_i$ and $D_{\Sigma}:=P\Sigma P^*$. Gathering these new notations yields 
$$\gamma_{n,p}(\theta')=\frac{1}{p^2}tr\left[\left(I_{p^2}-D(\theta')\right)D_{\Sigma}\overline{{\bf Y} {\bf Y}^*}\left(I_{p^2}-D(\theta')\right)\right]\ ,$$
where the vectors ${\bf Y}^i$ are independent standard Gaussian random vectors. Except $\overline{{\bf Y} {\bf Y}^*}$, every matrix involved in this last expression is diagonal. Besides, the diagonal matrix $D_{\Sigma}$ is positive since $\Sigma$ is non-singular. Thus, $$tr\left[(I_{p^2}-D(\theta'))D_{\Sigma}\overline{{\bf Y} {\bf Y}^*}(I_{p^2}-D(\theta'))\right]$$ is almost surely a positive quadratic form on the vector space generated by $I_{p^2}$ and $D(\Theta^+)$. Since the function $D(.)$ is injective and linear on $\Theta^+$, it follows that $\gamma_{n,p}(.)$ is almost surely strictly convex on $\Theta^+$.
\end{proof}\vspace{0.5cm}

\begin{proof}[Proof of Lemma 4.1 and Corollary 4.2 in \cite{verzelen_gmrf_theorie}]
The proof only uses the stationarity of the field $X$ on $\Lambda$ and the $l_1$ norm of $\theta$. However, the computations are a bit cumbersome. Let $\theta$ be an element of $\Theta^+$. By standard Gaussian properties, the expectation of $X{\scriptstyle[0,0]}$ given the remaining covariates is
\begin{eqnarray*}
  \mathbb{E}_{\theta}\left(X{\scriptstyle[0,0]}|X_{-\left\{0,0\right\}}\right) = \sum_{(i,j)\in \Lambda \backslash(0,0)}\theta{\scriptstyle[i,j]} X{\scriptstyle[i,j]}\ .
\end{eqnarray*}
By assumption $(\mathbb{H}_2)$, the $l_1$ norm of $\theta$ is smaller than one. We shall prove  by backward induction that for
any subset $A$ of $\Lambda\backslash\{(0,0)\}$ the matrix
$\theta^A$ uniquely defined by
$$ \mathbb{E}_{\theta}\left(X{\scriptstyle[0,0]}|X_A\right) = \sum_{(i,j)\in A} \theta^A{\scriptstyle[i,j]}X{\scriptstyle[i,j]}
\text{ and }\theta^A{\scriptstyle[i,j]}=0\text{ for any }(i,j)\notin A$$ satisfies $\|\theta^A\|_1\leq \|\theta\|_1$.
The property is clearly true if  $A=\Lambda\backslash\{(0,0)\}$. Suppose we have proved it for any set
 of cardinality $q$ larger than one. Let $A$ be a subset of $\Lambda\backslash\{(0,0)\}$ of cardinality
$q-1$ and $(i,j)$ be an element of $\Lambda\backslash (A\cup\{(0,0)\})$. 
Let us derive the expectation of $X{\scriptstyle[0,0]}$ conditionally to $X_A$ from the
expectation of $X{\scriptstyle[0,0]}$ conditionally to $X_{A\cup \{(i,j)\}}$.
\begin{eqnarray}
\mathbb{E}_{\theta}\left(X{\scriptstyle[0,0]}|X_A\right) & = & \mathbb{E}_{\theta}\left[\mathbb{E}(X{\scriptstyle[0,0]}|X_{A})|X_{A\cup
  \{(i,j)\}}\right] \nonumber \\ 
& =  & \sum_{(k,l)\in A}\theta^{A\cup\{(i,j)\}}{\scriptstyle[k,l]}X{\scriptstyle[k,l]} +
  \theta^{A\cup\{(i,j)\}}{\scriptstyle[i,j]}\mathbb{E}_{\theta}\left[X{\scriptstyle[i,j]}|X_A\right]\ . \label{decomposition_conditionelle1} 
\end{eqnarray}
Let us take the conditional expectation of $X{\scriptstyle[i,j]}$ with respect to 
$X_{A\cup\{(0,0)\}}$. Since the field $X$ is stationary on $\Lambda$ and by the induction hypothesis, the unique matrix
$\theta^{A\cup\{(0,0)\}}_{(i,j)}$ defined by
$$\mathbb{E}_{\theta}\left(X{\scriptstyle[i,j]}|X_{A\cup \{(0,0)\}}\right) = \sum_{(k,l) \in A\cup\{(0,0)\}}
\theta_{(i,j)}^{A\cup\{(0,0)\}}{\scriptstyle[k,l]} X{\scriptstyle[k,l]}$$  and
$\theta_{(i,j)}^{A\cup\{(0,0)\}}{\scriptstyle[k,l]}=0$ for any $(k,l)\notin A\cup\{(0,0)\}$ satisfies  $\|\theta_{(i,j)}^{A\cup\{(0,0)\}}\|_{1}\leq \|\theta\|_1$.
Taking the expectation conditionally to $X_A$ of this previous expression leads
to
\begin{eqnarray}
\mathbb{E}_{\theta}\left(X{\scriptstyle[i,j]}|X_{A}\right) = \sum_{(k,l) \in A}
\theta_{(i,j)}^{A\cup\{(0,0)\}}{\scriptstyle[k,l]} X{\scriptstyle[k,l]} +
\theta_{(i,j)}^{A\cup\{(0,0)\}}{\scriptstyle[0,0]}\mathbb{E}
\left(X{\scriptstyle[0,0]}|X_A\right)  \label{decomposition_conditionelle2}\ . 
\end{eqnarray}
Gathering identities (\ref{decomposition_conditionelle1}) and
(\ref{decomposition_conditionelle2}) yields 
\begin{eqnarray}
\mathbb{E}_{\theta}\left(X{\scriptstyle[0,0]}|X_A\right) = \sum_{(k,l)\in A}\frac{\theta^{A\cup\{i,j\}}{\scriptstyle[k,l]}
  +\theta^{A\cup\{(i,j)\}}{\scriptstyle[i,j]}\theta_{(i,j)}^{A\cup\{0,0\}}{\scriptstyle[k,l]}}
{1-\theta^{A\cup\{(i,j)\}}{\scriptstyle[i,j]}\theta_{(i,j)}^{A\cup\{0,0\}}{\scriptstyle[0,0]}}X{\scriptstyle[k,l]}\ ,  \nonumber
\end{eqnarray} 
since $\big|\theta^{A\cup\{i,j\}}{\scriptstyle[i,j]}\theta^{i,j}_{A\cup\{0,0\}}{\scriptstyle[0,0]}\big|< 1$. Then, we upper bound the $l_1$ norm of $\theta^{A}$ using that  
$\|\theta^{A\cup\{(i,j)\}}\|_1$ and $\|\theta^{A\cup\{(0,0)\}}_{(i,j)}\|_1$ are smaller or equal to $\|\theta\|_1$.
\begin{eqnarray*}
\lefteqn{\|\theta^{A}\|_1} &&\\& \leq & 
\frac{\sum_{(k,l)\in
  A}\left|\theta^{A\cup\{j+1\}}{\scriptstyle[k,l]}\right|+ \sum_{(k,l)\in A}
\left|\theta^{A\cup\{(i,j)\}}{\scriptstyle[i,j]}\theta_{(i,j)}^{A\cup\{0,0\}}{\scriptstyle[k,l]} \right|}{1-\left|\theta^{A\cup\{(i,j)\}}{\scriptstyle[i,j]}\theta_{(i,j)}^{A\cup\{(0,0)\}}{\scriptstyle[0,0]}\right|} \nonumber \\ 
& \leq &
\frac{\|\theta\|_1 +
\left|\theta^{A\cup\{(i,j)\}}{\scriptstyle[i,j]}\right|\left(\sum_{(k,l) \in
  A\cup\{(0,0)\}}\left|\theta^{A\cup\{(0,0)\}}_{(i,j)}{\scriptstyle[k,l]}\right|-1 - \left|\theta^{A\cup\{(0,0)\}}_{(i,j)}{\scriptstyle[0,0]}\right|
\right)}{1-\left|\theta^{A\cup\{(i,j)\}}{\scriptstyle[i,j]}\theta^{A\cup\{(0,0)\}}_{i,j}{\scriptstyle[0,0]}\right|} \nonumber\\ 
& \leq & \frac{\|\theta\|_1\left(1+ \left|\theta^{A\cup\{(i,j)\}}{\scriptstyle[i,j]}\right|\right) -
  \left|\theta^{A\cup\{(i,j)\}}{\scriptstyle[i,j]}\right|\left(1+ \left|\theta^{A\cup\{(0,0)\}}_{i,j}{\scriptstyle[0,0]}\right|\right)}{1-\left|\theta^{A\cup\{(i,j)\}}{\scriptstyle[i,j]}\theta^{A\cup\{(0,0)\}}_{(i,j)}{\scriptstyle[0,0]}\right|} \nonumber\\
& \leq & \|\theta\|_1 + \frac{\left|\theta^{A\cup\{(i,j)\}}{\scriptstyle[i,j]}\right|(\|\theta\|_1-1)\left(1+ \left|\theta^{A\cup\{(0,0)\}}_{(i,j)}{\scriptstyle[0,0]}\right|\right)}{1-\left|\theta^{A\cup\{(i,j)\}}{\scriptstyle[i,j]}\theta^{A\cup\{(0,0)\}}_{(i,j)}{\scriptstyle[0,0]}\right|}\ .
\end{eqnarray*}
Since $\|\theta\|_1$ is smaller than one, it follows that $\|\theta^A\|_1\leq \|\theta\|_1$.\\

Let $m$ be a model in the collection $\mathcal{M}_1$. Since $m$ stands for a set of neighbors of $(0,0)$, we may define $\theta^{m}$  as above. It follows that $\|\theta^{ m }\|_1\leq \|\theta\|_1$.
Since the field $X$ is stationary on the torus, $X$ follows the same distribution as the field $X^s$ defined by $X^s{\scriptstyle[i,j]}=X{\scriptstyle[-i,-j]}$. By uniqueness of
$\theta^{ m }$, we obtain that  $\theta^{ m }{\scriptstyle[i,j]}=\theta^{ m }{\scriptstyle[-i,-j]}$. Thus, $\theta^{ m }$ belongs to the space $\Theta_m$. Moreover, $\theta^{m}$ minimizes the function $\gamma(.)$ on $\Theta_m$. Since the $l_1$ norm of $\theta^m$ is smaller than one, $\theta^m$  belongs to $\Theta_{m,2}^+$. The matrices $\theta^{ m }$ and $\theta_{m,\rho_1}$ are therefore equal, which concludes the proof in the non-isotropic case.\\

Let us now turn to the isotropic case. Let $\theta$ belong to $\Theta^{\text{iso},+}$ and
let $m$ be a model in $\mathcal{M}_1$. As previously, the matrix $\theta^{ m }$ satisfies $\|\theta^{ m }\|_1\leq \|\theta\|_1$. Since the
  distribution of $X$ is invariant under the action of the group $G$,
  $\theta^{ m }$ belongs to $\Theta^{\text{iso}}_{m}$.  Since $\|\theta^{ m }\|_1\leq \|\theta\|_1$, $\theta^m$ lies in  $\Theta^{+,\text{iso}}_{m,2}$. It follows
  that $\theta^{ m } = \theta^{\text{iso}}_{m,\rho_1}$.
\end{proof}\vspace{0.5cm}

\begin{proof}[Proof of Corollary 4.3 in \cite{verzelen_gmrf_theorie}]
Let $\theta$ be a matrix in $\Theta^+$ such that  $(\mathbb{H}_2)$ holds and let $m$ be a model in $\mathcal{M}_1$. We decompose 
$\gamma(\widehat{\theta}_{m,\rho_1})$ using the conditional expectation of $X{\scriptstyle[0,0]}$ given $X_{ m }$.
\begin{eqnarray*}
\gamma(\widehat{\theta}_{m,\rho_1})&  =  & \mathbb{E}_{\theta}\bigg[X{\scriptstyle[0,0]} -
  \sum_{(i,j)\in  m }\widehat{\theta}_{m,\rho_1}{\scriptstyle[i,j]}X{\scriptstyle[i,j]}\bigg]^2\\
& = & \mathbb{E}_{\theta}\bigg[X{\scriptstyle[0,0]} -
  \mathbb{E}_{\theta}\left(X{\scriptstyle[0,0]}\left|X_{ m }\right.\right)\bigg]^2\\
&  +&	
  \mathbb{E}_{\theta}\bigg[\mathbb{E}_{\theta}\left(X{\scriptstyle[0,0]}\left|X_{ m }\right.\right)-
  \sum_{(i,j)\in  m }\widehat{\theta}_{m,\rho_1}{\scriptstyle[i,j]}X{\scriptstyle[i,j]} \bigg]^2. 
\end{eqnarray*}
By Corollary (11) in \cite{verzelen_gmrf_theorie}, we know that  $$\mathbb{E}_{\theta}\left(X{\scriptstyle[0,0]}\left|X_{ m }\right.\right)= \sum_{(i,j)\in  m }\theta_{m,\rho_1}{\scriptstyle[i,j]}X{\scriptstyle[i,j]}\ .$$ Combining these two last identities yields 
\begin{eqnarray*}
\gamma(\widehat{\theta}_{m,\rho_1})&  =  &  \gamma(\theta_{m,\rho_1}) +
\mathbb{E}_{\theta}\bigg[\sum_{(i,j)\in\Lambda\backslash\{(0,0)\}}\left(\theta_{m,\rho_1} -\widehat{\theta}_{m,\rho_1}\right){\scriptstyle[i,j]}X{\scriptstyle[i,j]}\bigg]^2\ .
\end{eqnarray*}	
Subtracting $\gamma(\theta)$, we obtain the first result. The proof is analogous
in the isotropic case.
\end{proof}

\appendix


\bibliographystyle{acmtrans-ims}

\bibliography{spatial}

\end{document}